\documentclass[titlepage,11pt,a4paper]{article}

\usepackage{amstext}
\usepackage{amssymb}
\usepackage{amsfonts}
\usepackage{amssymb}
\usepackage{amsmath}
\usepackage{graphicx}
\usepackage[arrow,matrix,curve]{xy}
\usepackage[latin1]{inputenc}
\usepackage{theorem}
\newtheorem{satz}{Theorem}[section]
\newtheorem{defi}[satz]{Definition}
\newtheorem{bem}[satz]{Remark}
\newtheorem{kor}[satz]{Corollary}
\newtheorem{lem}[satz]{Lemma}

\newtheorem{bei}[satz]{Example}

\newcommand{\dju}{\bigcup\hspace{-2.23mm\dot}\hspace{2mm}}
\newcommand{\modu}{~\rm{mod}~}
\newcommand{\smdju}{\cup\hspace{-2.15mm\dot}\hspace{2mm}}
\newcommand{\smdjuk}{\cup\hspace{-1.84mm\dot}\hspace{2mm}}

\newcommand{\qed}{\begin{flushright}$\square$\end{flushright}}
\newcommand{\filt}[6]{
\[
\begin{xy}
\xymatrix@R20pt@C20pt{
&\mathbb{C}^3&\\\langle #1,#2 \rangle\ar[ru]&\langle #3,#4\rangle\ar[u]&\langle #5,#6\rangle\ar[lu]\\\langle #1\rangle\ar[u]&\langle #3\rangle\ar[u]&\langle #5\rangle\ar[u]}
\end{xy}
\]
} 

\begin{document}
\newpage
\begin{center} 
{\bf\LARGE Tree modules of the generalized Kronecker quiver\\}
\end{center}
\vspace{0.0cm}
\begin{center}
Thorsten Weist\\Fachbereich C - Mathematik\\
Bergische Universität Wuppertal\\
D - 42097 Wuppertal, Germany\\
e-mail: weist@math.uni-wuppertal.de
\end{center}
\vspace{0.0cm}
\section{Introduction}
Fixing some representation of a quiver, we can choose a basis of the vector spaces associated to each vertex of the quiver and consider the maps restricted to these basis elements. We investigate the {\it coefficient quiver} in which the basis elements label the vertices and which has an arrow between two vertices if the matrix coefficient corresponding to these two basis elements is not zero. A representation is a called {\it tree module} if there exists a basis such that the coefficient quiver is a tree. This leads us to the following problem stated by Ringel, $\cite{rin2}$:\\
Does there exist an indecomposable tree module for every wild hereditary quiver and every root $d$? In particular, Ringel conjectured that there should be more than one isomorphism class for imaginary roots.\\
In this paper we prove the existence of indecomposable tree modules for the $m$-Kronecker quiver with $m\geq 3$, which is an extension of a result of $\cite{fahr}$, where the $3$-Kronecker quiver with dimension vector $(d,e)$ with $d<e<2d$ is treated. Moreover, we give an explicit construction of the coefficient quivers of indecomposable tree modules.\\\\
For coprime dimension vectors tree modules are obtained by considering the so called localization method in the Kronecker moduli space. This moduli space is a smooth projective variety parametrizing the stable representations and comes with an action of the $m$-dimensional torus. By results of $\cite{rei4}$ and $\cite{wei}$, there is a weight space decomposition of the vector spaces if a stable representation of the Kronecker quiver is a torus fixed point. This weight space decomposition provides the basis for the decomposition into basis vectors. Since this method also explicitly describes the maps between these weight spaces, and since the fixed points may be understood as representations of the universal abelian covering quiver, which is very similar to the regular $m$-tree, the coefficient quiver can be constructed from this.\\
Based on $\cite{wei}$, the construction of stable simple quivers, which are bipartite quivers with only two subquivers at the boundary, implies the existence of stable tree modules for coprime dimension vectors $(d,e)$ such that $d<e\leq (m-1)d+1$. Applying the reflection functor, see $\cite{bgp}$, we obtain all coprime cases. The simple quivers are a special case of the construction presented in the second section. It describes how to get new stable bipartite quivers by glueing two stable bipartite quivers of smaller dimension types where these quivers and dimension types have to satisfy certain conditions.\\
After proving the coprime case, the general case can be treated. In these cases we always get tree modules as factor modules of tree modules with coprime dimension vector.\\\\
{\bf Acknowledgment:} I would like to thank Philipp Fahr and Markus Reineke for very helpful discussions. Moreover, I would like to thank Claus Michael Ringel for his very useful comments.
\section{Torus fixed points in Kronecker moduli spaces}
After fixing notation for quivers and their representations, the fixed points of the Kronecker moduli space under a torus action are described. In particular, we study the stable simple quivers corresponding to dimension vectors $(d,e)$ of the Kronecker quiver with $d<e\leq (m-1)d+1$. Their stable representations correspond to torus fixed points and they will be shown to be tree modules. 
\subsection{Generalities}\label{allg}
Let $k$ be an algebraically closed field.
\begin{defi}
A quiver $Q$ consists of a set of vertices $Q_0$ and a set of arrows $Q_1$ denoted by $\alpha:i\rightarrow j$ for $i,j\in Q_0$.\\A vertex $q\in Q_0$ is called sink if there does not exist an arrow $\alpha:q\rightarrow q'\in Q_1$.\\
A vertex $q\in Q_0$ is called source if there does not exist an arrow $\alpha:q'\rightarrow q\in Q_1$.\\
A quiver is finite if $Q_0$ and $Q_1$ are finite.\\
A quiver is bipartite if $Q_0=I\smdjuk J$ such that for all arrows $\alpha:i\rightarrow j$ we have $i\in I$ and $j\in J$.
\end{defi}
Define the abelian group
\[\mathbb{Z}Q_0=\bigoplus_{i\in Q_0}\mathbb{Z}i\] and its monoid of dimension vectors $\mathbb{N}Q_0$.\\
A finite-dimensional $k$-representation of $Q$ is given by a tuple
\[X=((X_i)_{i\in Q_0},(X_{\alpha})_{\alpha\in Q_1}:X_i\rightarrow X_j)\]
of finite-dimensional $k$-vector spaces and $k$-linear maps between them. The dimension vector $\underline{\dim}X\in\mathbb{N}Q_0$ of $X$ is defined by
\[\underline{\dim}X=\sum_{i\in Q_0}\dim_kX_ii.\]
Let $d\in\mathbb{N}Q_0$ be a dimension vector. The variety $R_d(Q)$ of $k$-representations of $Q$ with dimension vector
$d$ is defined as the affine $k$-vector space
\[R_d(Q)=\bigoplus_{\alpha:i\rightarrow j} \mathrm{Hom}_k(k^{d_i},k^{d_j}).\]
The algebraic group 
\[G_d=\prod_{i\in I} Gl_{d_i}(k)\]
acts on $R_d(Q)$ via simultaneous base change, i.e.
\[(g_i)_{i\in Q_0}\ast
(X_{\alpha})_{\alpha\in Q_1}=(g_jX_{\alpha}g_i^{-1})_{\alpha:i\rightarrow
j}.\] The orbits are in bijection with the isomorphism classes of
$k$-representations of $Q$ with dimension vector $d$.\\\\
In the space of $\mathbb{Z}$-linear functions $\mathrm{Hom}_{\mathbb{Z}}(\mathbb{Z}Q_0,\mathbb{Z})$ we consider the basis given by the elements $i^{\ast}$ for $i\in Q_0$, i.e.
$i^{\ast}(j)=\delta_{i,j}$ for $j\in Q_0$. Define
\[\dim:=\sum_{i\in Q_0}i^{\ast}.\]
After choosing $\Theta\in
\mathrm{Hom}_{\mathbb{Z}}(\mathbb{Z}Q_0,\mathbb{Z})$,
we define the slope function
$\mu:\mathbb{N}Q_0\rightarrow\mathbb{Q}$ via
\[\mu(d)=\frac{\Theta{(d)}}{\dim(d)}.\]
The slope $\mu(\underline{\dim}X)$ of a representation $X$ of $Q$ is abbreviated to $\mu(X)$.
\begin{defi}
A representation $X$ of $Q$ is semistable (resp. stable) if for all subrepresentations  $U\subset X$ (resp. all proper subrepresentations $0\neq U\subsetneq X$) the following holds:
\[\mu(U)\leq\mu(X)\text{ (resp. } \mu(U)<\mu(X)).\]
\end{defi}
\subsection{The generalized Kronecker quiver}\label{KQ}
In the following let $k=\mathbb{C}$. The main focus of the paper is on the $m$-Kronecker quiver, $K(m)$, with $m\geq 3$. This is the quiver with two vertices and $m$ arrows between them:
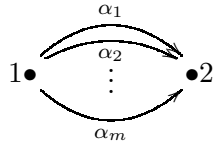
\[
\begin{xy}
\xymatrix@R10pt@C20pt{
\\1\bullet\ar@/_1.5pc/[rr]_{\alpha_m}\ar@/^1.5pc/[rr]^{\alpha_1}\ar@/^1.0pc/[rr]_{\alpha_2}&\vdots &\bullet 2\\&&
}
\end{xy}\]
A representation of $K(m)$ with dimension vector $(d,e)$ is given by $\mathbb{C}$-vector spaces $V$ and $W$ with dimensions $d$ and $e$ respectively and an $m$-tuple of $e\times d$-matrices 
\[(X_1,\ldots,X_m)\in\bigoplus_{i=1}^m \mathrm{Hom}(V,W)=:M_{e,d}(\mathbb{C})^m.\]\\ As described in the previous section, the group $G_{d,e}:=(Gl(V)\times Gl(W))$ acts on $M_{e,d}(\mathbb{C})^m$ via simultaneous base change. For $\Theta:=(1,0)$ the slope function
$\mu:\mathbb{N}^2\rightarrow\mathbb{Q}$ is given by
\[\mu(d,e)=\frac{d}{d+e}.\]
Thus we obtain the following (semi-)stability criterion for Kronecker representations:
\begin{lem}
A point $(X_1,\ldots,X_m)\in
M_{e,d}(\mathbb{C})^m$ is semistable if and only if for all subspaces  $U\subset V$ the following holds:
\begin{center}
$\dim\sum\limits_{k=1}^m X_k(U)\geq \dim U\cdot \frac{\textstyle
\dim W}{\textstyle \dim V}.$
\end{center}
It is stable if and only if for all proper subspaces
$0\neq U\subsetneq V$ the following holds:
\begin{center}
$\dim\sum\limits_{k=1}^m X_k(U)> \dim U\cdot \frac{\textstyle
\dim W}{\textstyle \dim V}.$
\end{center}
\end{lem}
Thus $d$ and $e$ being coprime implies that semistable points are stable. In this case we denote by $M_{d,e}^m$ the Kronecker moduli space, i.e. the smooth projective variety parametrizing the isomorphism classes of stable representations, for more details see $\cite{kin}$.\\Let $T:=(\mathbb{C}^{\ast})^{m}$ be the $m$-dimensional torus. It acts on $R_d(K(m))$ via
\[(t_1,\ldots,t_m)\cdot(X_1,\ldots,X_m)=(t_1X_1,\ldots,t_mX_m).\] 
Since the torus action commutes with the $G_{d,e}$-action, it induces a $T$-action on $M_{d,e}^m$. The torus fixed points can be described as follows, see $\cite{rei4}$ and $\cite{wei}$:
\begin{lem}
Let $X$ be a torus fixed point. Then there exist direct sum decompositions
\[V=\bigoplus_{\chi\in\mathbb{Z}^m}V_{\chi}\] and \[W=\bigoplus_{\chi\in
\mathbb{Z}^m}W_{\chi}\] 
such that
\[X_k(V_{\chi})\subseteq W_{\chi+e_k}\]
for all $\chi\in\mathbb{Z}^m$ and $k=1,\ldots ,m$.
\end{lem}
Here $e_1,\ldots,e_m$ denotes the standard basis of $\mathbb{Z}^m$.\\\\
The universal abelian covering quiver
$\hat{Q}$ has vertices $(1,\chi)$ and $(2,\chi)$, where $\chi$ runs through all $\chi\in\mathbb{Z}^m$, and arrows \[(1,\chi)\rightarrow (2,\chi+e_k)\] for each $k\in\{1,\ldots,m\}$  and $\chi\in\mathbb{Z}^m$.\\
For every fixed point there exists a unique dimension vector $\hat{d}$ of $\hat{Q}$ given by
\[\hat{d}_{1,\chi}=\dim V_{\chi}\text{ and
}\hat{d}_{2,\chi}=\dim W_{\chi}\] for $(1,\chi),(2,\chi)\in\hat{Q}_0$. Thus, by the preceding lemma, we may consider a fixed point of $M^m_{d,e}$ as a representations of $\hat{Q}$.
\begin{bem}
\end{bem}
\begin{itemize}
\item
The stability condition for representations of $\hat{Q}$ is induced from the original linear form $\Theta=(1,0)$. It is given by
\[\mu(\hat{d})=\frac{\sum_{\chi\in\mathbb{Z}^m}\hat{d}_{1,\chi}}{\sum_{\chi\in\mathbb{Z}^m}(\hat{d}_{1,\chi}+\hat{d}_{2,\chi})}.\]
If we consider bipartite quivers in general, we will also use the slope function induced by the linear form $\Theta=(1,0)$.
\end{itemize}
Consider the group action of $\mathbb{Z}^m$ on $\hat{Q}_0$ defined as follows:
\[\mu\cdot (i,\chi)=(i,\chi + \mu)\]
for $i=1,2$ and $\chi\in\mathbb{Z}^m$. This induces a group action on the set of dimension vectors $\mathbb{N}\hat{Q}_0$. Two dimension vectors in the same orbit of this action are said to be equivalent.\\\\
The other way around a stable representation of $\hat{Q}$ with dimension vector $\hat{d}$ corresponds to a torus fixed point with dimension vector $(d,e)$ where
\[d=\sum_{\chi\in\mathbb{Z}^m}\hat{d}_{1,\chi}\]
and
\[e=\sum_{\chi\in\mathbb{Z}^m}\hat{d}_{2,\chi}\]
respectively. A dimension vector $\hat{d}$ fulfilling these properties is called compatible with $(d,e)$. For a general bipartite quiver with a fixed dimension vector we will also refer to $(d,e)$ defined in this way as dimension type of a bipartite quiver.\\
In conclusion we have the following theorem, again see $\cite{rei4}$ and $\cite{wei}$:
\begin{satz}\label{isom}
The set of fixed points $(M_{d,e}^m)^T$ is isomorphic to the disjoint union of moduli spaces
\[\bigcup_{\hat{d}} M^s_{\hat{d}}(\hat{Q}),\] in which $\hat{d}$ ranges over all equivalence classes of dimension vectors being compatible with $(d,e)$.
\end{satz}
\subsection{Stability of bipartite quivers}\label{stabilitaet}
Now we will investigate stable quivers arising from the localization method presented in the last section. We will see that each colouring of the arrows which satisfies certain properties gives rise to a torus fixed point of the Kronecker quiver. In particular, we will study how to construct new stable quivers by glueing stable quivers of smaller dimension types. Moreover, we will treat stable simple quivers as a special case. Representations of such quivers are also torus fixed points and, moreover, indecomposable tree modules of the Kronecker quiver. Finally, we will construct stable tree modules for all roots except $(d,kd)$ and their reflections.
\begin{defi}
A tuple consisting of a quiver and a dimension vector is called stable if there exists at least one stable representation for this quiver and dimension vector.
\end{defi}
If it is clear which dimension vector we consider, we will simply call such a tuple stable quiver.
\begin{bem}\label{darstbip}
\end{bem}
\begin{itemize}
\item In order to test a $m$-bipartite quiver with
\[|\{\alpha\in Q_1\mid \alpha:i\rightarrow j\}|\leq 1\]
for stability, we do not need to consider an explicit representation. Fixing a dimension vector we can rather consider an arbitrary representation $X$ satisfying for all $j\in J$ and all subsets $A'_j\subseteq A_j$ with $R'_j:=|A'_j|$ the following property:
\[\dim(\bigcap_{i\in A'_j}X_{\alpha}(X_i))=\max\{0,\sum_{i\in A'_j}\dim(X_{\alpha}(X_i))-(R'_j-1)\dim(X_j)\}.\]
Indeed, if we consider a bipartite quiver of the form
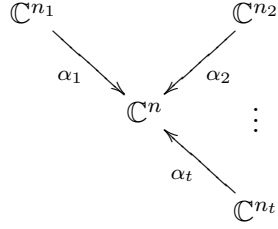
\[
\begin{xy}
\xymatrix@R20pt@C20pt{
\mathbb{C}^{n_1}\ar[rd]_{\alpha_1}&&\mathbb{C}^{n_2}\ar[ld]^{\alpha_2}\\
&\mathbb{C}^{n}&\vdots\\
&&\mathbb{C}^{n_t}\ar[lu]^{\alpha_t}}
\end{xy}
\]
with $n_i\leq n$ for all $1\leq i\leq t$, there always exists a representation of this quiver such that for all tuples of linear maps $X_{\alpha_{i_1}},\ldots,X_{\alpha_{i_k}}$ with $1\leq k\leq t$ and $1\leq i_1<i_2<\ldots<i_k\leq t$ the dimension of the intersection of their images is minimal. One verifies the existence and the dimension formula by induction on the number of arrows.
\end{itemize}
For a bipartite quiver $Q$ with vertex set $I\cup J$ and dimension vector $d$ define the sets
\[A_i:=\{j\in J\mid \alpha:i\rightarrow j\in Q_1,d_j\geq 1\}\]
and
\[A_j:=\{i\in I\mid \alpha:i\rightarrow j\in Q_1,d_i\geq 1\}.\]
Furthermore, define $R_i=|A_i|$ and $R_j=|A_j|$. 
A bipartite quiver is $m$-bipartite if we have for all sources $i\in I$ and all sinks $j\in J$ that
\[R_i,R_j\leq m.\]
Let $Q$ be a stable $m$-bipartite quiver of dimension type $(d,e)$ and arrow set $R$.
\begin{bem}\label{anzko}
\end{bem}
\begin{itemize} 
\item Let $i\in I$ be a vertex such that $\dim(i)=1$. Then we have
$m\geq R_i>\frac{e}{d}$.
\item For all $j\in J$ we have $R_j\leq m$.
\item Consider some colouring of the arrows $c:R\rightarrow\{1,\ldots,m\}$ satisfying the conditions:
For all $(i,j),(i,j')\in R$ such that $j\neq j'$ we have $c(i,j')\neq c(i,j)$ and for $(i,j),(i',j)\in
R$ we have $c(i,j)\neq c(i',j)$. In this way, we get a subquiver of the quiver $\hat{K}(m)$ defined in the last section. Therefore, by Theorem $\ref{isom}$ every stable representation of $Q$ defines a representation of the Kronecker quiver. We call a colouring satisfying these conditions stable.
\end{itemize}
Let $Q=(I\cup J,Q_1)$ and $Q'=(I'\cup J',Q_1')$ be two bipartite quivers with $j\in J$, $j'\in J'$. Define the bipartite quiver \[Q_{j,j'}(Q,Q')=(I\smdju I'\smdju J\backslash j \smdju J'\backslash j'\smdju j'' ,Q_1'')\] 
such that $\alpha:i\mapsto j_1\in Q_1''$ if and only if $\alpha:i\mapsto j_1\in Q_1$ or $\alpha:i\mapsto j_1\in Q_1'$ with $j_1\neq j,j'$ and  $\alpha:i\mapsto j''\in Q_1''$ if and only if  $\alpha:i\mapsto j_1\in Q_1$ or $\alpha:i\mapsto j_1\in Q_1'$ such that $j_1=j$ or $j_1=j'$.\\
Thus the new quiver is generated by the former ones by identifying two vertices of the set of sinks of these quivers.
\begin{defi}
The quiver $Q_{j,j'}(Q,Q')$ is called the glueing quiver of $Q$ and $Q'$ and the vertices $j,j'=j''$ the glueing vertices.
\end{defi}
\begin{defi}
Let $Q$ be a bipartite quiver with sources $I$. A subquiver of $Q$ with sources $I'$ is called boundary quiver if there exists precisely one $i_0\in I'$ such that $|A_{i_0}\cap A_{I\backslash I'}|=1$ and $|A_i\cap A_{I\backslash I'}|=0$ for all $i\neq i_0$. A boundary quiver is called proper boundary quiver if it does not contain any other boundary quiver.
\end{defi}
This means that boundary quivers are such subquivers which only have one common sink with the remainder of the quiver.\\In what follows we abbreviate the dimension of the image of a subspace $U$ to $d_U$. Thus if $U=\bigoplus_{i\in I}U_i$, we define
\[d_U:=\sum_{\alpha:i\rightarrow j}X_{\alpha}(U_i).\]
For a given dimension vector $(d,e)$ we now determine a unique dimension vector $(d_s,e_s)$ such that we are able to construct new stable quivers of dimension type $(d_s+kd,e_s+ke)$ by glueing quivers of the types $(d_s,e_s)$ and $k(d,e)$.\\
Fixing some dimension vector $(d,e)$, we first show that there exists a dimension vector $(d_s,e_s)$ such that $d_s\leq d$ and $e_s\leq e$ satisfying the conditions
\begin{itemize}
\item $\frac{e+e_s}{d+d_s}d<e+1$
\item $\frac{e+e_s}{d+d_s}d>e$
\item $\frac{e_s-1}{d_s}\leq \frac{e}{d}$ if $d\neq 1$ and $(e_s-1)d=ed_s$ if $d=1$ 
\item $\frac{e+e_s}{d+d_s}d'<\lceil\frac{e}{d}d'\rceil\hspace{0.3cm}\forall d'<d$
\item $\gcd(d+d_s,e+e_s)=1.$
\end{itemize}
We refer to this conditions as glueing condition. The first property is equivalent to the following:
\begin{eqnarray*}
 de+de_s &<& de+d+d_se+d_s\\
\Leftrightarrow  de_s &<& d+d_se+d_s\\
\Leftrightarrow  d(e_s-1)&<& d_s(e+1)\\
\Leftrightarrow  \frac{e_s-1}{d_s} &<&\frac{e+1}{d}.\\
\end{eqnarray*} 
The second one is equivalent to:
\begin{eqnarray*}
ed+e_sd &\geq ed+ed_s\\
\Leftrightarrow \frac{e_s}{d_s}>\frac{e}{d}.
\end{eqnarray*}
Therefore, it suffices to prove the second and third property because the first one follows from the third one.
\begin{lem}\label{startexist}
Let $(d,e)\in\mathbb{N}^2$ such that $d\leq e$ and $d,e$ are coprime. There exists a dimension vector $(d_s,e_s)$ satisfying the glueing condition. It is uniquely determined if we also assume $d_s\leq d$ and $e_s\leq e$. 
\end{lem}
{\it Proof.}
We first consider the special case $d=1$. It is easy to see that $(0,1)$ satisfies these properties for $(d,e)=(1,n)$ with $n\in\mathbb{N}$.\\\\
If $d\geq 2$, we already have $e\geq 3$. Choose $d_s\in\mathbb{N}$ minimal such that \[d\mid 1+ed_s. \]
This is possible because $\gcd(d,e)=1$ and, therefore, there exist $\lambda ',\mu '$ such that 
\[\lambda 'd=1-\mu 'e.\]
If $\mu '>0$, we have
\[ \lambda '^2d^2=1-2\mu 'e+\mu '^2e^2=1+e(\mu '^2e-2\mu ').\]
Since $\mu '^2e> 2\mu '$ for $e>2$, we obtain the existence and in particular that $d_s\in\mathbb{N}$.\\
Define
\[e_s=\frac{1+e(d+d_s)-de}{d}=\frac{1+d_se}{d}.\]
Because of the choice of $d_s$, we have $e_s\in\mathbb{N}$.\\Moreover, we get
\[-e(d+d_s)+d(e+e_s)=-ed-ed_s+de+d_se+1=1.\]
It follows that $\gcd(d+d_s,e+e_s)=1$.\\
Now we get
\[\frac{e_s}{d_s}=\frac{1+d_se}{dd_s}>\frac{e}{d}\]
and also
\[\frac{e_s-1}{d_s}=\frac{d_se-d+1}{dd_s}<\frac{e}{d}.\]
Thus it remains to prove the fourth property. By an easy calculation we get
\[\frac{e+e_s}{d+d_s}=\frac{e}{d}\left(\frac{ed+ed_s+1}{ed+ed_s}\right)=\frac{e}{d}\left( 1+\frac{1}{ed+ed_s}\right) .\]
Moreover, since 
\[\lceil\frac{e}{d}d'\rceil-\frac{e}{d}d'\geq\frac{1}{d}\]
 and 
\[\frac{d'}{ed+d_se}<\frac{1}{d+d_s}\]
for each $d'<d$, the existence of such a vector follows.
\begin{flushright}$\square$\end{flushright}
In what follows, we call such a vector $(d_s,e_s)$ satisfying these properties starting quiver for $(d,e)$. Now we show that $(d_s+kd,e_s+ke)$ with $k\geq 1$ also is a starting vector for the dimension vectors $(ld,le)$ with $l\geq 1$. Thus we can glue quivers of dimension type $(ld,le)$ on a stable quiver of type $(d_s+kd,e_s+ke)$ in order to obtain new stable quivers of dimension type $(d_s,e_s)+(k+l)(d,e)$.
\begin{kor}
Let $d,$ $d_s$, $e$ and $e_s$ fulfil the glueing condition. Then we have $\gcd(d_s+kd,e_s+ke)=1$ for all $k\geq 1$.
\end{kor}
{\it Proof.}
As before we also have 
\[-e(kd+d_s)+d(ke+e_s)=1\]
for arbitrary $k\geq 1$.
\begin{flushright}$\square$\end{flushright}
\begin{kor}\label{dvek}
Let $d,d_s,e,e_s$ satisfy the glueing condition and let $k,l\in\mathbb{N}$. We have 
\begin{enumerate}
\item \[\frac{e_s+ke}{d_s+kd}ld<le+1\]
\item \[\frac{e_s+ke}{d_s+kd}ld> le\]
\item \[\frac{e_s+ke-1}{d_s+kd}\leq\frac{le}{ld}\]
\item \[\frac{ke+e_s}{kd+d_s}d'<\lceil\frac{e}{d}d'\rceil\hspace{0.3cm}\forall d'<d\]
\end{enumerate}
\end{kor}
{\it Proof.}
The second (resp. third) statement is equivalent to the second (resp. third) property of the dimension vector $(d+d_s,e+e_s)$ what can be seen easily. Now the first property again follows from the third one. Let $k>1$. Then we have
\[\frac{e_s+ke}{d_s+kd}\leq\frac{e+e_s}{d+d_s}\]
which one verifies by an easy calculation. The fourth property follows from this.
\qed
\begin{bem}\label{darstellung}
\end{bem}
\begin{itemize}
\item If we want to decompose a dimension vector $(d,e)$ into
\[(d,e)=(d_s,e_s)+k(d',e')\]
for coprime $(d,e)$ such that $(d',e')$ and $(d_s,e_s)$ satisfy the glueing condition, we can proceed as follows:\\
fix $e'$ minimal such that
\[e\mid 1+de'\] and
\[d'=\frac{1+e'd}{e}.\] Now we compute $d_s$ and $e_s$ as before. It can be seen easily that these numbers satisfy the glueing condition. Indeed, one checks that
\[\frac{e-e_s}{e'}=\frac{d-d_s}{d'}.\]
From this it follows that $e'\mid e-e_s$ and $d'\mid d-d_s$ 
because $\gcd(d',e')=1$ and, trivially, $e-e_s,d-d_s\in\mathbb{N}$ hold. Now define $k=\frac{d-d_s}{d'}$. 
\end{itemize}
We need other properties of these natural numbers. By use of
\[e_sd-ed_s=1,\]for $0\leq k'\leq k$ we get
\[(ke+e_s)(k'd+d_s)+k-k'=(kd+d_s)(k'e+e_s).\]
For $d_1=k'd+d'\in\mathbb{N}$ with $0\leq d'< d$ and $0<d_1\leq kd+d_s$ define a map
\[f(d_1)=\min\{n\in\mathbb{N}\mid \frac{(ke+e_s)d_1+n}{kd+d_s}\in\mathbb{N}\}.\]
Note that $f$ is injective because $\gcd(d_s+kd,e_s+ke)=1$.
Then we get the following lemma:
\begin{lem}
Let $d_s$, $e_s$, $d$, $e$ fulfil the glueing condition. Then we have
\[(ke+e_s)(k'd+d_s)+k-k'=0\modu (kd+d_s)\]
for all $k'\leq k$.\\
Let $d_1=k'd+d'$ with $0\leq d'<d$. In particular, we have $f(d_1)=k-k'$ if $d'=d_s$ and thus $f(d_1)\geq k+1$ if $d'\neq d_s$. 
\end{lem}
Now we show how to get a stable quiver of dimension type $(d_s+(k+l)d,e_s+(k+l)e)$ by glueing a stable quiver of type $(d_s+kd,e_s+ke)$ and certain quivers of type $(ld,le+1)$. 
At this, Corollary $\ref{dvek}$ assures that the considered subquivers satisfy the stability condition. Afterwards we prove that also direct summands of representations do not contradict the stability condition what again follows by Corollary $\ref{dvek}$. We again point out Remark $\ref{darstbip}$. Thus we do not consider specific representations, but those satisfying the properties mentioned in the remark.\\\\
Fix natural numbers  $d$, $e$ and $m$ and let $\mathcal{S}_{ld,le+1}^m$ be the set of tuples consisting of a $m$-bipartite quiver of dimension type $(ld,le+1)$ and a sink $j$ satisfying the following properties:
\begin{itemize}
\item There exists at most one arrow between two vertices.
\item After decreasing the dimension of the sink $j$ by one, the resulting quiver is connected and semistable.
\item There exists a representation for the quiver such that for every $d'$-dimensional subspace $U$ we have
\[d_U>\frac{(k+l)e+e_s}{(k+l)d+d_s}d'.\]
\end{itemize}
Let $\mathcal{T}_{d,e}^m$ the set of all stable $m$-bipartite quivers of dimension type $(d,e)$.
\begin{satz}\label{verkl}
Let $d,d_s,e,e_s$ fulfil the glueing condition and let $k\in\mathbb{N}$. Let $T^0\in\mathcal{T}^m_{d_s+kd,e_s+kd}$ and $(T^1,j_1)\in\mathcal{S}_{ld,le+1}^m$. Moreover, let $j_0$ be a sink of $T_0^0$ such that $R_{j_0}+R_{j_1}\leq m$. Then $Q_{j_0,j_1}(T^0,T^1)$ with glueing vertex $j_2$ where $\dim(j_2):=\dim(j_0)+\dim(j_1)-1$ is an element of $\mathcal{T}^m_{d_s+(k+l)d,e_s+(k+l)e}$.
\end{satz}
{\it Proof.} 
For some subspace $U$ of one of the two subquivers we denote by $d_U$ the dimension of its image corresponding to its original quiver and by $d'_U$ the dimension of its image corresponding to the glueing quiver.\\
First let $U$ be a $d'$-dimensional subspace of $T^1$ such that $d'<ld$. Then by definition we have \[\frac{(k+l)e+e_s}{(k+l)d+d_s}d'<d_U=d'_U.\]
If $d'=ld$, the same inequality follows from $d_U=le+1$ together with the first property of Corollary $\ref{dvek}$.\\
Since we also have
\[\frac{e_s+ke}{d_s+kd}>\frac{e_s+(k+l)e}{d_s+(k+l)d},\] see the properties of the dimension vectors, the same follows for subspaces of the subquiver $T^0$.\\
It remains to prove that subspaces composed of subspaces of both subquivers fulfil the stability condition. Thus let $U'$ and $U''$ be two subspaces of dimension $1\leq d'\leq ld$ and $1\leq d''\leq kd+d_s$ such that we have proper inequality at least once.\\
Now it suffices to prove that
\[d'_{U'\oplus U''}\geq\frac{le}{ld}d'+d_{U''}>\frac{(k+l)e+e_s}{(k+l)d+d_s}(d'+d'')\]
where the first inequality follows from the semistability of the quiver obtained from $T^1$ after decreasing the dimension of the vertex $j_1$ by one. This is equivalent to
\[d_{U''}>\frac{(k+l)e+e_s}{(k+l)d+d_s}d''+\frac{d'}{d((k+l)d+d_s)}\]
using $e_sd-d_se=1$.\\
By the preceding lemma together with the assumption we have
\[d_{U''}\geq \frac{(ke+e_s)d''+f(d'')}{kd+d_s}.\]
First let $d''<kd+ds$. Assuming without lose of generality that $d'=ld$, it remains to prove that
\[ld''+((k+l)d+d_s)f(d'')>l(kd+d_s).\]
But this is easily verified.\\
Finally, let $d''=kd+d_s$ and $d'=l'd+d_1<ld$ with $0\leq d_1<d$.
We have
\[\frac{(k+l)e+e_s}{(k+l)d+d_s}(kd+d_s)=ke+e_s-\frac{l}{(k+l)d+d_s}\]
again using $e_sd-ed_s=1$. Thus it remains to prove 
\[\lceil\frac{e}{d}(l'd+d_1)\rceil=l'e+\lceil\frac{ed_1}{d}\rceil>\frac{(k+l)e+e_s}{(k+l)d+d_s}(l'd+d_1)-\frac{l}{(k+l)d+d_s}\]
what follows from the fourth property of Corollary $\ref{dvek}$ together with $l>l'$.
\qed
If $T^0$ and $T^1$ satisfy the condition of the theorem we call $T^0$ starting quiver for $T^1$.\\
Next, we apply the result to specific quivers. Therefore, let $T\in \mathcal{T}^m_{d,e}$. Starting with this quiver, we construct new quivers $\hat{T}$ of dimension type $(d,e+1)$ in one of the following ways:
\begin{itemize} 
\item Choose an $i\in I$ such that $R_i<m$ and define the new quiver by the vertex set $\hat{T}_0=T_0\cup\{ j\}$ and the arrow set $\hat{T}_1=T_1\cup\{\alpha:i\rightarrow j\}$. Finally, let $\dim(j)=1$.\item Choose a vertex $j\in J$ with $1<R_j<m$ and increase the dimension of the vertex by one.
\item Choose a vertex $j\in J$ such that \[\dim(j)<\sum_{i\in A_{j}} \dim(i)\] 
and increase the dimension of the vertex $j$ by one.
\end{itemize}
Denote the set of the resulting quivers by $\hat{\mathcal{T}}_{d,e}^m$ and refer to $j$ as modified vertex.
\begin{kor}\label{korverkl}
Let $d,d_s,e,e_s$ be as before and let $k\in\mathbb{N}$. Moreover, let $T^0\in\mathcal{T}^m_{d_s+kd,e_s+ke}$ and $T^1\in\hat{\mathcal{T}}^m_{d,e}$ with modified vertex $j_1$. Further let $j_0$ be a sink of $T_0^0$ such that $R_{j_0}+R_{j_1}\leq m$. Then $Q_{j_0,j_1}(T^0,T^1)$ with glueing vertex $j$, where $\dim(j):=\dim(j_0)+\dim(j_1)-1$, is an element of $\mathcal{T}_{d_s+(k+1)d,e_s+(k+1)e}$.
\end{kor}
{\it Proof.}
Let $U$ be a $d'$-dimensional subspace of $T^1$. Since $T^1$ results from a stable quiver we have $d_U>\frac{e}{d}d'$. Moreover, by the fourth property of Corollary $\ref{dvek}$ it follows that \[\frac{(k+1)e+e_s}{(k+1)d+d_s}d'<\lceil\frac{e}{d}d'\rceil\leq d_U\leq d'_U.\]
If $d'=d$, the same inequality follows from the first property together with \[\dim(j_1)\leq\sum_{i\in A_j} \dim(i)\]
and
\[d_U=e+1>\frac{(k+1)e+e_s}{(k+1)d+d_s}d.\]
\qed
Fixing a coprime dimension vector $(d,e)$ we now deal with the question how to construct a certain set of stable quivers. Therefore, we assign a set of stable quivers to tuple of natural numbers which is uniquely determined by the dimension vector, see also Example $\ref{tupel}$. These numbers correspond to the number of possible glueing vertices and possible colourings of the constructed quivers. They may be used to get a lower bound for the number of stable quivers and, therefore, for the Euler characteristic of Kronecker moduli spaces, see $\cite{wei}$. But, since we are only interested in the existence of tree modules, we will just apply this method to the special case of simple quivers.\\
Fix a dimension vector $(d,e)$ and the corresponding starting vector $(d_s,e_s)$. Denote by $\mathcal{T}^{(d,e)}_{n_1}$ the set of stable quivers of dimension type $(d_s,e_s)+n_1(d,e)$ with $n_1\geq 1$. As before let $\hat{\mathcal{T}}^{(d,e)}_{n_1}$ be the set which results by modifying a vertex $j_1$. Now we continue recursively: let $S\in\mathcal{T}^{(d,e)}_{n_k-1,\ldots,n_1}$ and $T\in\hat{\mathcal{T}}^{(d,e)}_{n_k,\ldots,n_1}$. Now let $\mathcal{T}^{(d,e)}_{1,n_k,\ldots,n_1}$ be the set consisting of all quivers $Q_{j_0,j_1}(S,T)$ such that $R_{j_0}+R_{j_1}\leq m$. Moreover, let the dimension of the glueing vertex $j$ be given by $\dim(j)=\dim(j_0)+\dim(j_1)-1$. In general let $\mathcal{T}^{(d,e)}_{n_{k+1},\ldots,n_1}$ 
be the set of glueing quivers resulting from glueing a quiver $S\in \mathcal{T}^{(d,e)}_{n_{k+1}-1,n_k,\ldots,n_1}$ and a quiver $T\in\hat{\mathcal{T}}^{(d,e)}_{n_k,\ldots,n_1}$ as described.
\begin{kor}\label{qst1}
The sets $\mathcal{T}^{(d,e)}_{n_k,\ldots,n_1}$ only contain stable quivers.
\end{kor}
{\it Proof.}
It suffices to prove that these quivers satisfy the condition of Corollary $\ref{korverkl}$.\\
We assume that $\mathcal{T}^{(d,e)}_{n_k,\ldots,n_1}$ only contains stable quivers. We have to prove that $\mathcal{T}^{(d,e)}_{n_{k+1},\ldots,n_1}$ just consists of stable quivers for all $n_{k+1}\geq 1$. Therefore we show that the quivers in $\mathcal{T}^{(d,e)}_{n_k-1,\ldots,n_1}$ are starting quivers for quivers in $\hat{\mathcal{T}}^{(d,e)}_{n_k,\ldots,n_1}$.\\
Let $(d^k,e^k)$ be the dimension vector corresponding to $\hat{\mathcal{T}}^{(d,e)}_{n_k,\ldots,n_1}$ and $(d_s^k,e_s^k)$ the one belonging to $\mathcal{T}^{(d,e)}_{n_k-1,\ldots,n_1}$. It suffices to prove that
\[(d^{k+1}_s,e^{k+1}_s)=(d_s^k,e_s^k)+(n_k-1)(d^k,e^k)\]
is the starting vector for
\[(d^{k+1},e^{k+1})=(d_s^{k},e_s^{k})+n_k(d^{k},e^{k}).\]
Indeed, the quivers in $\hat{\mathcal{T}}^{(d,e)}_{n_k,\ldots,n_1}$ are obtained by the modification described in Corollary $\ref{korverkl}$. But this is equivalent to
\[e_s^{k+1}=\frac{1+d_s^{k+1}e^{k+1}}{d^{k+1}}\]
with the additional condition $d^{k+1}_s\leq d^{k+1}$, see Lemma $\ref{startexist}$. The second property follows immediately, the first one is equivalent to \[e_s^k=\frac{1+d_s^ke^k}{d^k},\]
what follows by a direct calculation. Therefore, the claim follows by the induction hypothesis.
\qed
\begin{bei}\label{tupel}
\end{bei}
Let $(d_s,e_s)=(0,1)$ and $(d,e)=(1,n-1)$. Then we obtain the corresponding tuple of natural numbers $(n_k,\ldots,n_1)$ to a fixed dimension vector by proceeding as mentioned in Remark $\ref{darstellung}$. More detailed we have $(d^k,e^k)=(d_s,e_s)+n_k(d^{k-1},e^{k-1})$ and in this way we recursively obtain the whole tuple. The recursion terminates if $(d_s,e_s)=(0,1)$.\\
Consider for instance $(d,e)=(8,13)$. The tuple of numbers is given by $(n_3,n_2,n_1)=(1,2,2)$ with $n=2$. Thus we get
\begin{eqnarray*}(d,e)&=&(3,5)+(5,8)=(1,2)+(2,3)+((1,2)+2(2,3))\\
&=&(0,1)+(1,1)+(0,1)+2(1,1)\\
&&+((0,1)+(1,1)+2((0,1)+2(1,1))).
\end{eqnarray*}
Initially, consider the stable quivers of the dimension types $(1,2)$ and $(2,3)$, i.e.
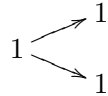
\[
\begin{xy}
\xymatrix@R0.5pt@C20pt{
&1\\1\ar[ru]\ar[rd]&\\&1}
\end{xy}
\]
and
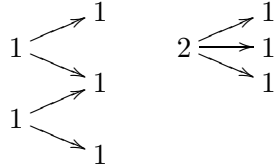
\[
\begin{xy}
\xymatrix@R0.5pt@C20pt{
&1&&1\\1\ar[ru]\ar[rd]&&2\ar[ru]\ar[r]\ar[rd]&1\\&1&&1\\1\ar[ru]\ar[rd]&\\&1 }
\end{xy}
\]
By use of Corollary $\ref{korverkl}$ we obtain the following stable quivers of dimension type $(3,5)$ by glueing:
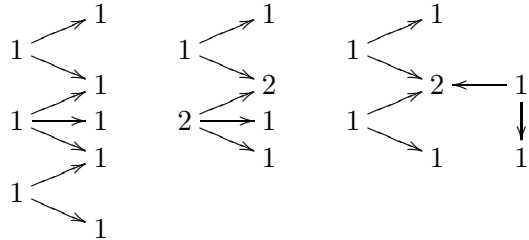
\[
\begin{xy}
\xymatrix@R0.5pt@C20pt{
&1&&1&&1\\1\ar[ru]\ar[rd]&&1\ar[ru]\ar[rd]&&1\ar[ru]\ar[rd]\\&1&&2&&2&1\ar[l]\ar[dd]\\1\ar[ru]\ar[rd]\ar[r]&1&2\ar[r]\ar[ru]\ar[rd]&1&1\ar[ru]\ar[rd]&\\&1&&1&&1&1\\1\ar[ru]\ar[rd]\\&1 }
\end{xy}
\]
Obviously, all modules of theses quivers are tree modules. Next, for instance we obtain the following stable quivers of type $(5,8)$ by glueing:
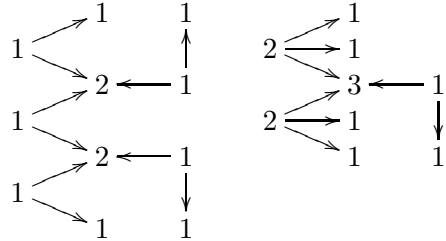
\[
\begin{xy}
\xymatrix@R0.5pt@C20pt{
&1&1&&1\\1\ar[ru]\ar[rd]&&&2\ar[ru]\ar[r]\ar[rd]&1\\&2&1\ar[l]\ar[uu]&&3&1\ar[dd]\ar[l]\\1\ar[ru]\ar[rd]&&&2\ar[ru]\ar[r]\ar[rd]&1\\&2&1\ar[dd]\ar[l]&&1&1\\1\ar[ru]\ar[rd]\\&1&1 }
\end{xy}
\]
Modules of the first quiver are trees. But for modules of the latter one this is not true in general. Finally, for instance the following bipartite quiver of dimension type $(8,13)$ is stable by use of Corollary \ref{korverkl}. Moreover, all modules are tree modules.
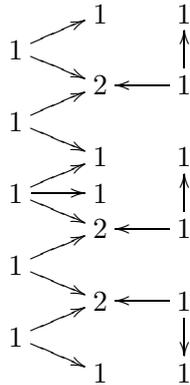
\[
\begin{xy}
\xymatrix@R0.5pt@C20pt{
&1&1\\1\ar[ru]\ar[rd]\\&2&1\ar[l]\ar[uu]\\1\ar[ru]\ar[rd]\\&1&1\\1\ar[ru]\ar[r]\ar[rd]&1\\&2&1\ar[uu]\ar[l]\\1\ar[ru]\ar[rd]\\&2&1\ar[l]\ar[dd]\\1\ar[ru]\ar[rd]\\&1&1
}
\end{xy}
\]

\subsection{Simple quivers}
In this section we investigate the class of so called simple quivers introduced in $\cite{wei}$. They are bipartite quivers with only two proper boundary quivers. The previous sections show that representations of these quivers are torus fixed points of the Kronecker quiver after some appropriate colouring of the arrows. It should be noted that these quivers are far away from being all bipartite quivers resulting from the localization method and the glueing method presented in the last section, see Example $\ref{tupel}$. Nevertheless, they give rise to an exponential growing class of stable representations which are also trees. But since the existence of the simple quivers is sufficient to prove the existence of tree modules for all dimension vectors, and the existence of more than one isomorphism class of them in the case of imaginary roots, we restrict to studying these quivers.\\
Finally, we show that there even exists a stable bipartite quiver for each dimension vector $(d,e)$ with $e\neq kd$ for $1\leq k\leq m$ and $d<e\leq (m-1)d+1$ whose representations are tree modules.\\\\
We assume that all vertices (before and after glueing) have dimension one, where in general the dimension of a vertex $q$ of a quiver $Q$ with a fixed dimension vector $d\in\mathbb{N}^{Q_0}$ is defined by the natural number $d_q$.\\Let $(Q^{l_k})_{k\in\mathbb{N}}$ be defined by \[Q^{l_k}_0=\{i_k\}\cup J_k\text{ and } Q^{l_k}_1=\{(i_k,j)\mid j\in J_k\}\] with $|J_k|=l_k$. Consider the glueing quiver \[Q^{l_1,l_2}:=Q_{j_1,j_2}(Q^{l_1},Q^{l_2})\] where $j_1\in J_1$ and $j_2\in J_2$. Also define the set of possible glueing vertices of the new quiver by $K_2:=J_2\backslash j_2$, where the glueing vertex is denoted by $j_2$.\\We continue recursively as follows: let $Q^{l_1,l_2,\ldots,l_n}$ and $Q^{l_{n+1}}$ be two bipartite quivers. Then we define 
\[Q^{l_1,\ldots,l_{n+1}}:=Q_{j_n,j_{n+1}}(Q^{l_1,\ldots,l_n},Q^{l_{n+1}})\] 
with $j_n\in K_n$ and $j_{n+1}\in Q^{l_{n+1}}_0$. Again define the set $K_{n+1}:=J_{n+1}\backslash j_{n+1}$.
\\\\Fix some $m\in\mathbb{N}$ such that $m\geq 3$. Let $n,t\in\mathbb{N}$ such that $2\leq n\leq m-1$. To each $(t+1)$-tuple $(s_1,\ldots,s_{t+1})_{m,n}\in\mathbb{N}^{t+1}$ we can assign the quiver \[Q^{n^{s_1},n+1,n^{s_2},n+1,\ldots,n+1,n^{s_{t+1}}},\] in which \[n^{s_i}:=\underbrace{n,\ldots,n}_{s_i-times}.\] 
Considering the case $m=3$ and $n=2$, we obtain quivers of the following form:\\
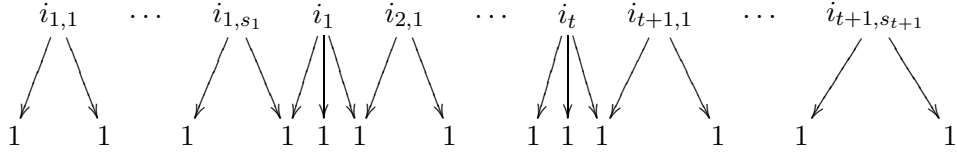
\[
\begin{xy}
\xymatrix@R30pt@C0.5pt{
&i_{1,1}\ar[ld]\ar[rd]&&\dots&&i_{1,s_1}\ar[ld]\ar[rd]&&i_1\ar[ld]\ar[rd]\ar[d]&&i_{2,1}\ar[ld]\ar[rd]&&\dots&&i_t\ar[ld]\ar[rd]\ar[d]&&i_{t+1,1}\ar[ld]\ar[rd]&&\dots&&i_{t+1,s_{t+1}}\ar[ld]\ar[rd]\\
1&&1&&1&&1&1&1&&1&&1&1&1&&1&&1&&1
}
\end{xy}
\]
where $\dim(i_k)=\dim(i_{j,k})=1$.
Similarly, define the quiver
\[\hat{Q}_{s}=Q^{n+1,n^{s_1},n+1,n^{s_2},\ldots,n+1,n^{s_{t+1}}}.\]
and
\[\hat{s}:=(s_1-1,s_2,\ldots,s_{t+1})_{m,n}.\]
Informally, we obtain the quiver $\hat{Q}_{\hat{s}}$ from the quiver $Q_s$ by glueing an extra arrow on ''the source at the left margin''. In the following we write $\hat{s}$ instead of $\hat{Q}_{\hat{s}}$ and $s$ instead of $Q_s$ if no confusion arises.
\begin{defi}
Let $m,n$ and a tuple $(s_1,s_2,\ldots,s_{t+1})_{m,n}$ be given. The corresponding quiver is called simple.
\end{defi}
We decompose the vertex set of $Q$ into $I\cup J$ such that the vertices in $I$ are the sources and the vertices in $J$ are the sinks.\\
Let $Q_1\subset I\times J$ be the set of arrows and let $X$ be a representation of some stable simple quiver. Let $c:Q_1\rightarrow\{1,\ldots,m\}$ be a stable $m$-colouring of the set of arrows.\\Then we obtain a representation of the $m$-Kronecker quiver in the obvious way. We set
\[V=\bigoplus_{i\in I}\mathbb{C}v_i\text{ and } W=\bigoplus_{j\in J}\mathbb{C}w_j\]
and $X_{\alpha_k}v_i=X_{\alpha}v_i=w_j$ if and only if there exists an arrow $\alpha:i\rightarrow j$ with $c(\alpha)=k$, where $k\in\{1,\ldots,m\}$.
The dimension vector of the corresponding Kronecker representation is given by
\[ d=\sum_{i=1}^{t+1} s_i + t\]
and
\[ e=(\sum_{i=1}^{t+1} (n-1)s_i+1) + (n-1)t.\]
From Section $\ref{KQ}$ and in particular by Theorem $\ref{isom}$ we get the following:
\begin{lem}
Let $s$ be a stable simple quiver with dimension vector $(d,e)$ and some stable colouring $c$. Then every stable representation $X$ of the simple quiver is a torus fixed point of the moduli space of stable representations via the construction described above. In particular $X$ is a stable representation of the Kronecker quiver.
\end{lem}
We test the simple quivers for stability, i.e. which of these quivers allow stable representations. Therefore, we can assume that $X_{\alpha}=1$ for every arrow $\alpha\in Q_1$.\\
Let $U$ be a subspace of $X$ with $\dim U=\sum_{i=k}^l s_i +(l-k)$. Since $X$ is stable, the dimension of the image satisfies the inequality
\[d_U\geq (\sum_{i=k}^l (n-1)s_i+1) + (n-1)(l-k).\]
It is easy to see that it is enough to test stability for those subspaces $U$ such that
\[d_U= (\sum_{i=k}^l (n-1)s_i+1) + (n-1)(l-k).\]
In this situation we have the following lemma:
\begin{lem}\label{stabs}
A quiver $(s_1,\ldots,s_{t+1})_{m,n}$ is stable if and only if 
\[d(l-k+1)>(t+1)(\sum_{i=k}^l s_i+(l-k))\]
for all $1\leq k< l\leq t+1$ with $l-k<t$.
\end{lem}
{\it Proof.} We have $e=(n-1)d+t+1$. Therefore, from the stability condition we obtain
\[(\sum_{i=k}^l (n-1)s_i+1) + (n-1)(l-k)>\frac{(n-1)d+t+1}{d}(\sum_{i=k}^l s_i +(l-k)).\]
Thus we get
\begin{eqnarray*}
&&d(\sum_{i=k}^l (n-1)s_i)+d(n-1)(l-k)+d(l-k+1)\\&>&(n-1)d(\sum_{i=k}^l s_i +(l-k))+(t+1)(\sum_{i=k}^l s_i +(l-k))\end{eqnarray*}
and therefore the assertion follows.
\begin{flushright}$\square$\end{flushright}
Define \[s_{k,l}:=(\sum_{i=k}^l s_i+l-k)(t+1).\]
Obviously, the stability condition is independent of $m$ and $n$. Thus define
\[(s_1,\ldots,s_{t+1}):=\{(s_1,\ldots,s_{t+1})_{m,n}\mid n\leq m\}.\]
Define $l_{(d,e)}:=\lceil \frac{e}{d}l\rceil$. It is easy to see that for each $l$-dimensional subspaces corresponding to a boundary quiver of a simple quiver we have $d_U=l_{(d,e)}$. Indeed, otherwise we would get a contradiction to
\[\lceil\frac{e}{d}l\rceil+\lceil\frac{e}{d}(d-l)\rceil-1=e.\]
\begin{lem}
Let $X$ be a stable representation of a simple quiver. For all subspaces $U=\oplus_{i\in I'} X_i$ of dimension $l$ with $I'\subset I$ we either have $d_U=l_{(d,e)}$ or $d_U=l_{(d,e)}+1$.
\end{lem}
{\it Proof.} Let $I_1,I_2\subset I$. We assume that there exist two subspaces $U_1=\oplus_{i\in I_1} X_i$ and $U_2=\oplus_{i\in I_2} X_i$ of the same dimension such that $d_{U_1}=d_{U_2}+t$ where $t\geq 2$. We may without lose of generality assume that $d_{U_2}=l_{(d,e)}$. Let $I_3:=I\backslash I_1$ and consider the subspace \[U_3=\bigoplus_{i\in I_3}X_i.\] Then we obtain $d_{U_3}=e-l_{(d,e)}.$\\Because of the stability condition we also get that
\[ (d-l)_{(d,e)}\leq e-l_{(d,e)}.\] It is checked by an easy calculation that this is only possible if $\frac{e}{d}l\in\mathbb{N}$. But this is impossible because of $\mathrm{gcd}(d,e)=1$ and $l<d$.
\begin{flushright}$\square$\end{flushright}
We get the following corollaries:
\begin{kor}\label{kk+1}
For a stable simple quiver $(s_1,\ldots,s_{t+1})$ there exists a $k\in\mathbb{N}$ such that for all $1\leq i\leq t+1$ we either have $s_i=k$ or $s_i=k+1$.
\end{kor}
{\it Proof.} Without lose of generality let $m=3$ and $n=2$. If $s_i=k$ for some $i$, there exists a subspace $U$ of dimension $k+2$ such that $d_U=k+5$. If there existed some $j$, such that $s_j=k+2$, there would also exist a subspace $U'$ of dimension $k+2$ such that $d_{U'}=k+3$. But such a subspace cannot exists because of the preceding lemma.
\qed
\begin{kor}\label{01}
Let $(s_1,\ldots,s_{t+1})$ be stable and let $k\in\mathbb{Z}$ such that $k\geq -\min_{i=1}^{t+1} s_i$. Then $(s_1+k,\ldots,s_{t+1}+k)$ is also stable. In particular, we may assume that $s_i=0$ or $s_i=1$.
\end{kor}
{\it Proof.}
By adding $k$ to each $s_i$ on both sides, we obtain the equivalence to the original condition.
\begin{flushright}$\square$\end{flushright}
\begin{kor}\label{min}
Let $(s_1,\ldots,s_t)$ be stable and
\[\sum_{i=1}^k s_i=a\] for some $a\in\mathbb{N}$ such that $k<t$. Then we have
\[\sum_{i=l}^{l+k-1} s_i\leq a\] for $l+k-1\leq t$.
\end{kor}
\begin{kor}
Let $(s_1,\ldots,s_t)$ be stable. Then we have $s_i=s_{t-i+1}$ for all $i\leq\frac{t}{2}$.
\end{kor}
{\it Proof.}
It is easy to see that the stability of $(s_1,\ldots,s_t)$ is equivalent to the stability of $(s_t,\ldots,s_1)$.\\
Now the assertion follows from the preceding corollary. Indeed, if there existed a non-symmetric simple quiver, it would follow that \[\sum_{i=1}^k s_i=a_1\text{, but also that } \sum_{i=t-k+1}^t s_i=a_2\]
where $a_1\neq a_2$. But this is a contradiction.
\begin{flushright}$\square$\end{flushright}
Now we can simplify the stability condition for simple quivers:
\begin{lem}\label{stabs2}
Let $s=(s_1,\ldots,s_{t+1})$ be a symmetric simple quiver.
Then the following are equivalent:
\begin{enumerate}
\item The quiver $s=(s_1,\ldots,s_{t+1})$ is stable.  
\item We have \[dl>s_{1,l}\] for all $1\leq l < t+1$. 
\item We have \[dl>s_{1,l}>dl-(t+1)\] for all $1\leq l \leq\frac{t+1}{2}$ 
\end{enumerate}
\end{lem}
{\it Proof.}
Assume that the third condition holds. Let $k,l$ such that $k\leq l$.
Because of the symmetry of $s$ we have
\[s_{t-k+2,t+1}=s_{1,k}.\]
Moreover, we have
\[s_{1,t+1}=d(t+1).\]
Thus it follows
\begin{eqnarray*} s_{k,l} &=& s_{1,t+1}-s_{1,k-1}-s_{l+1,t+1}-2(t+1)\\
&=& s_{1,t+1}-s_{1,k-1}-s_{1,t-l+1}-2(t+1)\\
&<& d(t+1)-d(k-1)-d(t-l)\\
&=& d(l-k+1).
\end{eqnarray*}
Now we show the equivalence of the second and third condition.
Thus assume that we have \[dl>s_{1,l}\] for all $1\leq l < t+1$. 
Using the symmetry we have for $k\leq\frac{t+1}{2}$ that
\begin{eqnarray*}
s_{1,k} &=& s_{1,t+1}-s_{1,t-k+1}-(t+1)\\
&>& d(t+1)-d(t-k+1)-(t+1)\\
&=& dk-(t+1).
\end{eqnarray*} 
The other way around let $\frac{t+1}{2}\leq l<t+1$. Since $t+1-l\leq\frac{t+1}{2}$, we have
\begin{eqnarray*}
s_{1,l} &=& s_{1,t+1}-s_{l+1,t+1}-(t+1)\\
&=& s_{1,t+1}-s_{1,t+1-l}-(t+1)\\
&<& d(t+1)-d(t+1-l)\\
&=& dl.
\end{eqnarray*}
This completes the proof.
\begin{flushright}$\square$\end{flushright}
The next theorem shows the existence and uniqueness of simple quivers for coprime dimension types.
\begin{satz}\label{eindeutigkeit}
Let $\mathrm{gcd}(d,e)=1$. Then there exists exactly one stable quiver $s\in\mathbb{N}^t$ of dimension type $(d,e)$.
\end{satz}
{\it Proof.}
Initially, we show the existence of such a quiver. Because of Corollary $\ref{01}$ we can without lose of generality assume that $t+1<d<2(t+1)$. Indeed, if $d=t$, we only get $s=0$. For the quiver $s$ we thus have
\[dl-(t+1)<s_{1,l}<dl\]
for all $l\leq\frac{t+1}{2}$.\\Since $\mathrm{gcd}(d,e)=1$ and $e=d+t+1$, it follows that $t+1$ is no factor of $d$. Therefore, such a quiver $s$ exists. Note that $t$ is already uniquely determined by $t=e-d$.\\
Now the uniqueness follows immediately. If $s'$ were another stable simple quiver of dimension type $(d,e)$, there would exist a minimal $k$ such that $s_k\neq s'_k$. In particular, we would have $s_k=s'_k+n$ for some $n\neq 0$. But then it would follow that either $s'_{1,k}>dk$ or $s'_{1,k}<dk-(t+1)$.
\begin{flushright}$\square$\end{flushright}
Now by the preceding theorem together with Corollary $\ref{qst1}$ we obtain the following:
\begin{kor}\label{satz1}
Let $d,d_s,e,e_s$ be natural numbers satisfying the glueing condition. Moreover, let $s_{(d_s,e_s)}$ be the stable simple quiver of dimension type $(d_s,e_s)$ and $s_{(d,e)}$ the one of dimension type $(d,e)$. The stable simple quiver of dimension type $(d_s+kd,e_s+ke)$ is given by $s_{(d_s+kd,e_s+ke)}=(s_{(d_s,e_s)},\hat{s}_{(d,e)}^k)$.
\end{kor}
In what follows, we assume that $s_i\in\{0,1\}$. This is possible because of Corollary \ref{01} and Corollary \ref{kk+1}.\\
Now we show how to construct the stable simple quivers recursively. Let $(n_1,\ldots,n_k)\in\mathbb{N}^k$ be a $k$-tuple of natural numbers. Define the following simple quivers: 
\begin{itemize}
\item $s_{n_1}=1^{n_1+1}$
\item $s_{n_1,n_2}=(s_{n_1-1},\hat{s}_{n_1}^{n_2})$
\item $s_{n_1,\ldots,n_{k+1}}=(s_{n_1,\ldots,n_k-1},\hat{s}_{n_1,\ldots,n_k}^{n_{k+1}})$
\end{itemize}
Denote in the following $s_{n_k}:=s_{n_1,\ldots,n_{k}}$ if no confusion arises.
Note that for $k\geq 2$ we obviously have $s_{n_1,\ldots,n_k,0}=s_{n_1,\ldots,n_{k}-1}$.\\
By Corollary $\ref{qst1}$ we obtain the following result:
\begin{kor}\label{qst2}
The quivers $s_{n_k}$ are stable.
\end{kor}
Denote
\[l^k:=\underbrace{l,\ldots ,l}_{k-times}\]
and consider the maps $\eta_{n}^l:\{l-1,l\}\rightarrow\mathbb{N}^n\cup\mathbb{N}^{n+1}$ and $\Theta_{n}^l:\{l-1,l\}\rightarrow\mathbb{N}^{n+1}\cup\mathbb{N}^{n+2}$ defined by
\begin{eqnarray*} \eta_{n}^l&:&(l-1)\mapsto (l-1),l^{n-1}\\
&&l\mapsto (l-1),l^n\\
\end{eqnarray*}
and
\begin{eqnarray*} \Theta_{n}^l&:&(l-1)\mapsto (l-1)^{n+1},l\\
&&l\mapsto (l-1)^n,l.\\
\end{eqnarray*}
Obviously we have $\Theta_n^l(l)=(\eta_1^l)^n(l)$. These maps are to be applied componentwise to vectors consisting of natural numbers $l-1$ and $l$ respectively. Denote $\eta_n=\eta_n^1$ and $\Theta_n=\Theta_n^1$ respectively.\\
We get the following lemma:
\begin{lem}
For all $k\in\mathbb{N}^+$ we have
\[\hat{s}_{n_k}=\eta_{n_1}\circ\eta_{n_2}\circ\ldots\circ\eta_{n_k}(1).\]
\end{lem}
{\it Proof.}
Let $k=1$. For an arbitrary $n_1$ we have $s_{n_1}=1^{n_1+1}$, thus $\hat{s}_{n_1}=01^{n_1}=\eta_{n_1}(1)$.\\Consider 
\[\hat{s}_{n_{k+1}}=(\hat{s}_{n_k-1},(\hat{s}_{n_k})^{n_{k+1}}).\]
Because of the induction hypothesis we obtain
\[\hat{s}_{n_k}=\eta_{n_1}\circ\eta_{n_2}\circ\ldots\circ\eta_{n_k}(1)\]
and, moreover, we have
\[\eta_{n_{k+1}}(1)=01^{n_{k+1}}.\]
Therefore, it suffices to prove
\[\hat{s}_{n_k-1}=\eta_{n_1}\circ\eta_{n_2}\circ\ldots\circ\eta_{n_k}(0).\]
For $k=1$ we have $\hat{s}_{n_1-1}=01^{n_1-1}=\eta_{n_1}(0)$. As before we get
\[\eta_{n_{k+1}-1}(0)=01^{n_{k+1}-1}.\]
Then it follows
\begin{eqnarray*}\hat{s}_{n_{k+1}-1}&=&(\hat{s}_{n_k-1},(\hat{s}_{n_k})^{n_{k+1}-1})\\
&=&(\eta_{n_1}\circ\eta_{n_2}\circ\ldots\circ\eta_{n_k}(0),(\eta_{n_1}\circ\eta_{n_2}\circ\ldots\circ\eta_{n_k}(1))^{n_{k+1}-1})\\
&=&(\eta_{n_1}\circ\eta_{n_2}\circ\ldots\circ\eta_{n_k}(01^{n_{k+1}-1}))\\
&=&\eta_{n_1}\circ\eta_{n_2}\circ\ldots\circ\eta_{n_{k+1}}(0).
\end{eqnarray*}
\begin{flushright}$\square$\end{flushright}
\begin{defi}
Let $s_{d,e}$ be the stable simple quiver of dimension type $(d,e)$ and let
\[\hat{s}_{d,e}=\eta_k\circ\ldots\circ\eta_1(l)\]
with $l\in\mathbb{N}$. The maps $\eta_1,\ldots,\eta_k$ are called quiver functions (of type $l$) for $d,e$.
\end{defi}
Now it suffices to determine the quiver functions (including its type) for a fixed dimension vector in order to get the unique stable simple quiver for a given coprime dimension vector.\\
We recursively obtain the quiver functions as follows: fix $d,e,m$ and $n$ with $(n-1)d<e<nd$ and $2\leq n\leq m-1$.
In this situation we have
\[(d,e)=k_{1,1}(1,n-1)+k_{1,2}(1,n).\]
After solving this system of linear equations, we have
\[k_{1,1}=nd-e\text{ and } k_{1,2}=e-(n-1)d.\]
If $k_{1,1}=0$, we are done, since in this case $k_{1,2}=1$ and consequently $\hat{s}_{d,e}=0$. Analogously, $k_{1,2}=0$ implies that $\hat{s}_{d,e}=1$.\\Therefore, assume $k_{1,2}\neq 0$. Then the type of the quiver functions is given by $l_1:=\lceil\frac{k_{1,1}}{k_{1,2}}\rceil$. If $k_{1,2}=1$, it follows that the type of the quiver functions is $k_{1,1}$ and we obtain $\hat{s}_{d,e}=k_{1,1}$. If $k_{1,1}=1$, we have $\hat{s}_{d,e}=0^{k_{1,2}-1}1$.\\
Thus assume $k_{1,2}\neq 1$. We recursively proceed as follows:\\if $k_{1,2}|k_{1,1}$, then $\gcd(d,e)=1$ implies that $k_{1,2}=1$. Indeed, otherwise $k_{1,2}$ would divide both $d$ and $e$. By considering the type of the quiver function, we get linear equations
\[(d,e)=k_{2,1}((l_1-1)(1,n-1)+(1,n))+k_{2,2}(l_1(1,n-1)+(1,n)).\]
From this it follows that
\[k_{2,1}=(l_1+1)e-(l_1n-l_1+n)d\]
and
\[k_{2,2}=d(l_1n-l_1+1)-l_1e.\]Define $l_2:=\lceil\frac{k_{2,2}}{k_{2,1}}\rceil$.
If $k_{2,2}=1$, we obtain $\hat{s}_{d,e}=\Theta_{k_{2,1}}(l)$. 
If $k_{2,2}\neq 1$, the first quiver function is of the form
\[\eta_{l_2}:l\mapsto (l-1)l^{l_2}.\]
If $k_{2,1}=1$, we have $\hat{s}_{d,e}=\eta_{l_2}(l)$.\\
If $k_{2,1},k_{2,2}\neq 1$, we proceed recursively by solving the systems of linear equations given by
\[k_{j,1}+k_{j,2}=k_{j-1,1}\]
and
\[(l_{j-1}-1)k_{j,1}+l_{j-1}k_{j,2}=k_{j-1,2}\]
as long as $k_{j,1}=1$ or $k_{j,2}=1$ where $l_{j-1}:=\lceil\frac{k_{j-1,2}}{k_{j-1,1}}\rceil$ and $j\geq 3$.
From this we obtain $\hat{s}_{d,e}$ as follows: if $k_{j,2}=1$, we have
\[\hat{s}_{d,e}=\eta_{l_2}\circ\ldots\circ\eta_{l_{j-1}}\circ\Theta_{k_{j,1}}(l)\]
and if $k_{j,1}=1$, we have
\[\hat{s}_{d,e}=\eta_{l_{2}}\circ\eta_{l_3}\circ\ldots\circ\eta_{l_{j}}(l).\]
Now we obtain from Corollary $\ref{qst2}$ that the resulting quivers $s_{d,e}^m$ are stable. Moreover, they are of the requested dimension type because the tuple $(k_{j,1},k_{j,2})$ exactly indicates the number of $(l-1)'s$ and $l's$ of the vector $\eta_{l_j}\circ\ldots \circ\eta_{l_k}(l)$.
\begin{bem}\label{tupelbest}
\end{bem}
\begin{itemize}
\item Consider Corollary $\ref{qst1}$ with $(d_s,e_s)=(0,1)$ and $(d,e)=(1,n-1)$ for some $n\geq 2$. Moreover, let $m\geq 3$ and fix a dimension vector $(d',e')$ such that $(m-1)d'\geq e'\geq d'$ and assume that $(n-1)d'\leq e'\leq nd'$ for some $n\in\mathbb{N}$.
By means of Corollary $\ref{satz1}$ it is checked that
\[(n_k,\ldots,n_1)=(l_k,\ldots,l_2,l+1).\]
Therefore, the tuple corresponding to a dimension vector is already determined by a system of linear equations. 
\end{itemize}
Finally, for every root $(d,e)\neq (d,kd)$ with $k\in\mathbb{N}$ satisfying $d<e\leq (m-1)d+1$ we construct a stable quiver which also is a tree module.\\
For a quiver $Q$ we define
\[N_i=\{j\in J\mid (\alpha:i\rightarrow j)\in Q_1\}\] for $i\in I$. We define $N_j$ analogously.\\
Fix $n\in\mathbb{N}$ such that $n<m$ and a tuple $s\in\mathbb{N}^{t+1}$ and consider the simple quiver corresponding to $(s_1,0^{t-1},s_{t+1})_{n,m}$. This quiver has a subquiver $Q^{(n+1)^{t-1}}$, which has exactly $(t-1)$ sources $q_1,\ldots,q_{t-1}$. Furthermore, for each source $q_i$ choose a sink $p_i\in N_{q_i}$ satisfying $|N_{p_i}|=1$. Now we can inductively construct the following quiver:
\[Q^{s_1,s_2,\ldots,s_k,s_{t+1}}=Q_{p_k,q}(Q^{s_1,s_2,\ldots,s_{k-1},s_{t+1}},Q^{n^{s_k}})\] where $2\leq k\leq t$ and $q\in Q^{n^{s_k}}_0$ is a sink of one of the at most two proper boundary quivers. In this way we can assign a quiver to each tuple $(s_1,\ldots,s_{t+1})\in\mathbb{N}^{t+1}$.
This quiver has the following stability condition:
\begin{lem}
Let $s=\sum_{i=1}^{t+1}s_i$. The quiver $(s_1,s_2,\ldots,s_{t+1})_{n,m}$ is stable if and only if
\[sl+t>t\sum_{i=1}^ls_i+l>s(l-1)+1\]
for all $l=1,\ldots,t$.
\end{lem}
{\it Proof.}
We have $(d,e)=(t-1+s,n(t-1)+(n-1)s+1)$. Fix $1\neq k\leq l\neq t+1$. Analogously to the proof of Lemma $\ref{stabs}$, by an easy calculation we obtain the following stability condition
\[s(l-k+2)+t>t\sum_{i=k}^ls_i+(l-k+2).\]
If $k=1$ and $l<t+1$, we get the condition
\[sl+t>t\sum_{i=1}^ls_i+l\]
and analogously for $k\neq 1$ and $l=t+1$ the condition
\[s(t+2-k)+t>t\sum_{i=k}^{t+1}s_i+(t+2-k).\]
From 
\[t\sum_{i=1}^ls_i+l=ts-t\sum_{i=l+1}^{t+1}s_i+l\]
the assertion follows similar to the proof of $\ref{stabs2}$.
\qed
\begin{bem}
\end{bem}
\begin{itemize}
\item It is easy to see that $(d,e)=(d,kd)$ is equivalent to $s=1$. Thus we do not get a stable quiver in this way.
\end{itemize}
Since $sl+t-s(l-1)-2=s+t-2\geq t$ for $s>1$, we obtain the following result in the same way as in the proof of Theorem $\ref{eindeutigkeit}$. But note that these quivers are not unique in general.
\begin{satz}\label{stabilekoecher}
For each dimension vector $(d,e)\neq (d,kd)$ with $d<e\leq (m-1)d+1$ there exists a tuple $(s_1,s_2,\ldots,s_{t+1})\in\mathbb{N}^{t+1}$ such that the corresponding quiver is stable.
\end{satz}
\section{Tree modules of the Kronecker quiver}
In this section indecomposable tree modules are constructed by means of the construction of stable torus fixed points in the last section. Stability of the corresponding representation of $\hat{Q}$ implies indecomposability of the representation, and the decomposition into weight spaces gives rise to a suitable basis. Moreover, an explicit method for constructing tree modules for all dimension vectors, by use of the reflection functor, is described.
\subsection{Coefficient quivers and tree modules}
\setcounter{subsection}{1}
In this section we introduce coefficient quivers and tree modules, following the presentation given in $\cite{rin1}$.\\
Let $Q$ be a quiver with dimension vector $d=(d_q)_{q\in Q_0}$ and let $X$ be a representation of $Q$. A basis of $X$ is a subset $\mathcal{B}$ of $\bigoplus_{q\in Q_0}X_q$ such that
\[\mathcal{B}_q:=\mathcal{B}\cap X_q\] is a basis of $X_q$ for all vertices $q\in Q_0$. For every arrow $\alpha:i\rightarrow j$ we may write $X_{\alpha}$ as a $(d_j\times d_i)$-matrix $X_{\alpha,\mathcal{B}}$ with coefficients in $\mathbb{C}$ such that the rows and columns are indexed by $\mathcal{B}_j$ and $\mathcal{B}_i$ respectively. If
\[X_{\alpha}(b)=\sum_{b'\in\mathcal{B}_j}\lambda_{b',b}b'\]
with $\lambda_{b',b}\in\mathbb{C}$, we obviously have $(X_{\alpha,\mathcal{B}})_{b',b}=\lambda_{b',b}$.
\begin{defi}
The coefficient quiver $\Gamma(X,\mathcal{B})$ of a representation $X$ with a fixed basis $\mathcal{B}$ has vertex set $\mathcal{B}$ and arrows between vertices are defined by the condition: if $(X_{\alpha,\mathcal{B}})_{b,b'}\neq 0$, there exists an arrow $(\alpha,b,b'):b\mapsto b'$.\\
A representation $X$ is called a tree module if there exists a basis $\mathcal{B}$ for $X$ such that the corresponding coefficient quiver is a tree.
\end{defi}
As already mentioned in the introduction, defining the coefficient quiver immediately raises the question posed by Claus Michael Ringel, $\cite{rin2}$:\\
Does there exist an indecomposable tree module for every wild hereditary quiver and every root $d$?\\He conjectured that there should be more than one isomorphism class for imaginary roots.
\subsection{Reflection of stable tree modules}
In this section we will apply the reflection functor to tree modules of the Kronecker quiver. We will see that reflected modules of tree modules are tree modules. Therefore, it suffices to prove the existence of tree modules for roots $(d,e)$ satisfying $d\leq e<(m-1)d$. Indeed every module has a corresponding module in this domain.\\\\
As usual denote by $E_q$ the simple representation corresponding to the vertex $q$ defined by $X_q=\mathbb{C}$ and $X_{q'}=0$ for all $q'\neq q$.\\For a quiver $Q$ consider the matrix $A=(a_{i,j})_{i,j\in Q_0}$ with $a_{i,i}=2$ and $-a_{i,j}=-a_{j,i}$ for $i\neq j$, in which $a_{i,j}=|\{\alpha\in Q_1\mid\alpha:i\rightarrow j\vee \alpha:j\rightarrow i\}|$.\\Fixed some $q\in Q_0$ define $r_q:\mathbb{Z}Q_0\rightarrow\mathbb{Z}Q_0$ as
\[r_q(q')=q'-a_{q,q'}q.\]
We have the following theorem, see $\cite{bgp}$ and $\cite{kac}$:
\begin{satz}
Let $Q$ be a quiver and $q\in Q_0$ a fixed vertex. Let $q$ be a sink (resp. a source). Then there exists a functor \[R_q^+ (resp. R_q^-):\modu\mathbb{C}Q\rightarrow \modu\mathbb{C}Q_q\] with the following properties (if $q$ is a source, replace $+$ by $-$):
\begin{enumerate}
\item $R_q^+(U\oplus U')=R_q^+(U)\oplus R_q^+(U')$\\
\item Let $U$ be an indecomposable representation of $Q$.
\begin{enumerate}
\item If $U\cong E_q$, then $R_q^+(E_q)=0$.\\
\item If $U\ncong E_q$, then $R_q^+(U)$ is indecomposable with $R_q^+R_q^-(U)\cong U$ and $\dim R_q^+(U)=r_q(\dim(U))$.
\end{enumerate}
\end{enumerate}
In particular, we have: $\mathrm{End} U\cong \mathrm{End} R_q^+(U)$.
\end{satz} 
We consider the infinite regular $m$-tree as a quiver and denote it by $T(m)$. Therefore, we fix a vertex and define it to be a source. Its adjacent vertices are defined to be sinks and we proceed recursively in this manner. In particular, every vertex is either a source or sink and not both.\\More formally:
Let $k_2\in\{1,2,\ldots,m\}$ and $k_l\in\{1,2,\ldots,m-1\}$ for some $l\in\mathbb{N}$ with $l\geq 3$. Define the index set $K_m$ as follows:
\[K_m=\{k=(1,k_2,k_3,\ldots,k_n)\mid n\geq 1\}.\]
Define $s(k)=n$ for $k=(1,k_2,k_3,\ldots,k_n)$ and $K_m^n=\{k\in K_m\mid s(k)=n\}$.\\
Consider the quiver $T(m)$ with vertex set $T(m)_0=I\dju J$, where $I=\{i_k\mid k\in K_m,s(k)\text{ odd}\}$ and $J=\{j_k\mid k\in K_m,s(k)\text{ even}\}$ whose set of arrows is given by: \[T(m)_1=\{\alpha:i\rightarrow j\mid i\in I\cap K_m^n,j\in J\cap K_m^{n+1}\}.\]\\
Let
\[X=((X)_{q\in T(m)_0},(X_{\alpha}:X_i\rightarrow X_j)_{\alpha\in T(m)_1})\]
be a representation. Then the dimension vector $\dim X$ is defined as \[\dim X=\sum_{q\in T(m)_0} \dim X_qq.\]\\
Let the slope function $\mu:\mathbb{N}T(m)_0\rightarrow\mathbb{Q}$ be given by
\[\mu(X)=\frac{\sum_{i\in I} \dim(X_i)}{\sum_{q\in T(m)_0} \dim(X_q)}.\]
For a representation $X$ of the regular $m$-tree define the integers \[d:=\sum_{i\in I}\dim(X_i)\text{ and } e:=\sum_{j\in J} \dim(X_j).\]
Obviously every representation of $T(m)$ can be considered as a representation of the quiver $\hat{Q}$ defined in $\ref{KQ}$ if we choose a stable colouring of the arrows of $T(m)$. Thus from \ref{isom} we get the following lemma:
\begin{lem}
Let $X$ be a stable representation of $T(m)$ with a stable colouring. Then $X$ corresponds to a torus fixed point of the Kronecker moduli space $M_{d,e}^m$, i.e. in particular some stable representation of $K(m)$ with dimension vector $(d,e)$.
\end{lem}
If we apply the reflection functor to some representation of the Kronecker quiver which is not isomorphic to $E_1$ and has dimension vector $(d,e)$, we obtain a representation with dimension vector $(e,me-d)$. Thereby we must bear in mind the fact that the arrows have to be turned around.
\begin{bem}\label{color}
\end{bem}
\begin{itemize}
\item
Assume $e>(m-1)d$. Considering the regular $m$-tree, every vertex is either sink or source of exactly $m$ arrows. In order to get a torus fixed point out of some stable representation of the regular $m$-tree these arrows have to be coloured such that the colours of the arrows are pairwise disjoint. It may happen that several vertices have the same weight. 
\item In the case $e<(m-1)d$ we consider stable simple quivers, which may obviously be coloured such that there does not exist any cycle. In particular, the resulting torus fixed points are tree modules.
\item
It suffices to test stability for indecomposable submodules. This is because if $U=U_1\oplus\ldots\oplus U_n$ is some submodule of some module $M$ such that every $U_i$ is indecomposable, then $\mu(U)\leq \max(\mu(U_1),\ldots,\mu(U_n))$. For more details see for instance $\cite{rei4}$.
\item Due to $\cite{kac}$ we get 
\[\frac{m-\sqrt{m^2-4}}{2}\leq\frac{e}{d}\leq\frac{m+\sqrt{m^2-4}}{2}\]
for every imaginary root $(d,e)$ of the Kronecker quiver. Due to the reflection functor every imaginary root $(d,e)$ has some corresponding dimension vector $(d',e')$ such that $d'<e'<(m-1)d'$.
\end{itemize}
Let $X$ be a representation of $T(m)$. Then define \[I_X=\{i\in I\mid \dim(X_i)\neq 0\vee (\exists j\in N_i:\dim(X_j)\neq 0)\}.\] This is a finite subset of $I$ whose cardinality $|I_X|$ we denote by $n_X$. Index the set $I_X=\{i_1,\ldots,i_{n_X}\}.$\\
Define $R_{I_X}^-=R^-_{i_1}\circ\ldots\circ R^-_{i_{n_X}}$.
Then the following holds:
\begin{lem}
Let $X$ be a stable representation of the regular $m$-tree $T(m)$ such that $X\ncong E_i$ for all $i\in I$. Then the representation $R_{I_X}^-(X)$ is also a stable representation of $T(m)$.
\end{lem}
{\it Proof.}
It is easy to see that the conclusion of the lemma is independent of the choice of the indexing of the set $I_X$.\\
Thus let $V$ be a submodule of $R_{I_X}^-(X)$, which may be assumed to be indecomposable due to the above remark. Moreover, we can assume that the submodule is of the form $(V_i,X_{\alpha}(V_i))$. These are exactly those submodules that are of the form $R_{I_X}^-(U)$ for some submodule $U$ of $X$. Indeed, we have:
\[\sum_{\alpha:i\rightarrow j\in A_i}R_q^-(X_{\alpha})(U_j)=\bigoplus_{\alpha:i\rightarrow j\in A_i}U_j/\mathrm{im}(h'|_{h'^{-1}(\bigoplus_{\alpha:i\rightarrow j\in A_i}U_j)}).\]
Consider the following two cases:\\
First let $V\cong E_i$ for some $i\in I$. Then we have $\mu(V)=0$ by the definition of stability for $Q$. Note that the arrows were turned around. But we also have $\mu(R_{I_X}^-(X))\neq 0$, otherwise we would have had $e=0$ (note again that the arrows were turned around) and therefore $X$ would be isomorphic to a direct sum of copies of $E_i$'s with $i\in I$. But this is contradictory to the assumption.\\
Now let $V\ncong E_i$ for all $i\in I$. Therefore, without lose of generality $V$ is of the form $R_{I_X}^-(U)$ for some submodule $U$ of $X$. Since $X$ is stable, we have $\mu(U)<\mu(X)$. Let $(d',e')$ be the dimension vector of $U$. Then the inequality implies $\frac{d'}{e'} < \frac{d}{e}$. Out of this one deduces by direct calculation that $\frac{e'}{me'-d'} < \frac{e}{me-d},$ but this means precisely that $\mu(V)<\mu(R_{I_X}^-(X))$.
\qed
\begin{bem}
\end{bem}
\begin{itemize}
\item
One easily checks that given a representation of the regular $m$-tree with dimension vector $(d,e)$, which corresponds to some representation of $K(m)$, that after applying the functor $R_{I_X}^-$ one gets in fact a representation of $K(m)$ with dimension vector $(e,me-d)$. Since all arrows are turned around (except such arrows that have zero-dimensional sink AND source after applying the functor), we have for the reflected dimension vector:
\begin{eqnarray*}r_{I_X}(d,e)&=&(\sum_{j\in J} \dim(X_j)j,m\sum_{j\in J}\dim(X_j)j-\sum_{i\in I}\dim(X_i)i)\\&=&(e,me-d),\end{eqnarray*}
because every vertex has exactly $m$ adjacent vertices.
\end{itemize}
If we just want to prove the existence of indecomposable tree modules it is enough to apply the following result, see \cite{rin1}:
\begin{satz}\label{ringel}
Let $Q$ be a quiver and let $X$ be a representation of $Q$ such that $\rm{Ext}_{\mathbb{C}Q}^1(X,X)=0$. Then $X$ is a tree module.
\end{satz}
In summary we get the following result:
\begin{satz}\label{teilerfremdBaum}
For every coprime dimension vector $(d,e)$ of the Kronecker quiver there exists an indecomposable tree module.
\end{satz}
{\it Proof.}
The case $e<d$ is covered by reversing arrows in the corresponding case $(e,d)$.\\
For $d<e\leq (m-1)d+1$ we get an indecomposable tree module by considering the stable simple quiver associated to this dimension vector and by applying an appropriate colouring, see $\ref{color}$. It is easy to see that this representation is exceptional, i.e. its \rm{Ext}-group vanishes. Indeed, we have 
\[\dim\mathrm{Hom}(X,X)-\dim\mathrm{Ext}(X,X)=\langle \dim X,\dim X\rangle =1\]
and $\dim\mathrm{Hom}(X,X)=1$.\\\\ 
Considering the case $e>(m-1)d$, we first determine the equivalent case such that $e'<(m-1)d'$. Then we apply the reflection functor suitably many times to some previously constructed module of the regular $m$-tree corresponding to a tree module of the Kronecker quiver. The reflected module also has an \rm{Ext}-group which vanishes. Thus by $\ref{ringel}$ it is also a tree module and we obtain the explicit coefficient quiver as described in Section $\ref{constr}$.
\qed
\subsection{The general case}
The case of arbitrary dimension vectors $d<e<2d$ for the $3$-Kronecker quiver is also explicitly described in $\cite{fahr}$. There an algorithm is described which states how to get tree shaped coefficient quivers for these dimension vectors starting from the dimension vector $(d,2d+1)$.\\
Consider the map $r:\mathbb{N}^2\rightarrow\mathbb{N}^2$ defined by
$r(d,e):=(e,me-d).$ We have the following result:
\begin{satz}
\begin{enumerate}
\item For every dimension vector $(d,e)$ there exists an indecomposable tree module of the Kronecker quiver.\item Let $k,l,n\in\mathbb{N}_0$. For each imaginary root $(d,e)\neq r^l(n,kn)$ there exists a stable tree module.
\end{enumerate}
\end{satz}
{\it Proof.}
Analogous to $\cite{fahr}$ we obtain from the stable simple quiver $s_{d,e}$ with $\gcd(d,e)=1$ all dimension vectors such that $d'\leq e'\leq (m-1)d'$ by removing arrows. To be precise: for $(d',e')$ consider the dimension vector $(d',e)$ with $e$ being minimal such that $\gcd(d',e)=1$ and $e\geq e'$. Obviously we have $(m-1)d'+1\geq e> d'$ because $\gcd(d',(m-1)d'+1)=1$. The stable simple quiver  $s_{d',e}$ is bipartite, thus we may decompose the vertex set in sources and sinks, denoted by $I\cup J$. Furthermore, we have $e-d'+1$ vertices $j\in J$ with $|N_j|=1$. Denote this set by $J'$. Choose from $J'$ exactly $e-e'$ vertices and denote this set by $J''$.\\
Define the quiver $s_{d',e'}$ as the quiver with vertex set $I\cup J\backslash J''$ and an arrow $\alpha:i\rightarrow j\in (s_{d',e'})_1$ if and only if $\alpha:i\rightarrow j\in (s_{d',e})_1$ and $j\in J\backslash J''$. Since there exists at least one stable (and in particular indecomposable) representation of $s_{d',e}$, there exists an indecomposable representation for $s_{d',e'}$ as well. This is because the representation remains indecomposable after removing these kinds of arrows. Note also that all vertices of this kind of stable quivers are supposed to be one-dimensional and that we get in this way a factor representation of $s_{d',e}$.\\As shown before we obtain all coprime cases via the reflection functor. The reflection functor maps this indecomposable factor representation again to an indecomposable factor representation, now of $R^-_{I_X}(s_{d',e})$. This factor representation is forced to be a tree because $R^-_{I_X}(s_{d',e})$ is one. Indeed, if there were a basis vector which were mapped to itself by a chain of maps, this chain would also exist in $R^-_{I_X}(s_{d',e})$.\\By this method we get all dimension vectors. Indeed, fixing some dimension vector $(d,e)$, we can determine the corresponding one such that $d\leq e\leq (m-1)d+1$ via the reflection $(d,e)\mapsto (md-e,d)$.
\\\\
We obtain the second part as follows. By Theorem \ref{stabilekoecher} for each dimension vector $(d,e)\neq (n,kn)$ with $d<e\leq (m-1)d+1$ there exists a stable tree module. 
The representations of this quiver have a vanishing $\rm{Ext}$-group which follows in the same way as in the proof of Theorem $\ref{teilerfremdBaum}$. Thus their reflected modules are tree modules as well.
\qed
\subsection{Construction of tree modules}\label{constr}
In order to get a concrete coefficient quiver, we first investigate how to apply the reflection functor to a constructed coefficient quiver.\\
Let $X$ be a stable representation of $T(m)$. Consider a sink $q$ with basis $\mathcal{B}_q$ and for each $\alpha:q_i\rightarrow q$ consider a basis  $\mathcal{B}_{q_i}$ with $1\leq i\leq m$. If there exists some $w\in\mathcal{B}_{q}$ such that $(X_{\alpha:q_k\rightarrow q,\mathcal{B}})_{v,w}=0$ for all $v\in\mathcal{B}_{q_k}$, we have $R_{q_k}^-(X_{\alpha:q_k\rightarrow q})(w)=\bar{w}$ with $\bar{w}\notin\mathcal{B}_{q_l}$ for all $1\leq l\leq m$. It is easy to see that we also have $(R_{q_k}^-(X_{\alpha:q_k\rightarrow q,\mathcal{B}}))_{w',\bar{w}}=0$ for each $w'\in\mathcal{B}_q$ with $w'\neq w$.
Therefore the determination of the coefficient quiver reduces to the cases in which $n$ arrows point from vertices of dimension one at one vertex of dimension $(n-1)$, i.e.:
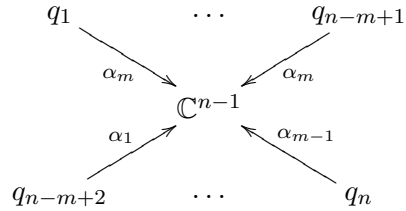
\[
\begin{xy}
\xymatrix@R20pt@C20pt{
q_1\ar[rd]_{\alpha_m}&\ldots&q_{n-m+1}\ar[ld]^{\alpha_m}\\
&\mathbb{C}^{n-1}&\\
q_{n-m+2}\ar[ru]^{\alpha_1}&\ldots&q_n\ar[lu]_{\alpha_{m-1}}}
\end{xy}
\]
where $q_i=\mathbb{C}$ for all $1\leq i\leq n$. This corresponds to the real root of the $n$-subspace quiver. By the indecomposability we have a coefficient quiver which looks as follows, see $\cite{fahr}$ for the case $m=3$ and $n=2$:
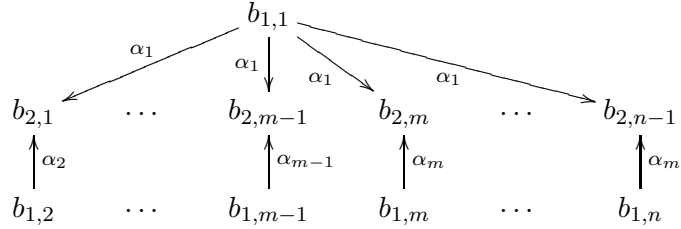
\[
\begin{xy}
\xymatrix@R20pt@C20pt{
&&b_{1,1}\ar[lld]_{\alpha_1}\ar[d]_{\alpha_1}\ar[rd]_{\alpha_1}\ar[rrrd]_{\alpha_1}&&&\\
b_{2,1}&\ldots&b_{2,m-1}&b_{2,m}&\ldots&b_{2,n-1}\\
b_{1,2}\ar[u]_{\alpha_2}&\ldots&b_{1,m-1}\ar[u]_{\alpha_{m-1}}&b_{1,m}\ar[u]_{\alpha_m}&\ldots&b_{1,n}\ar[u]_{\alpha_m}
}
\end{xy}
\]
Now the reflection functor behaves as mentioned above.
\begin{bem}\label{colo}
\end{bem}
\begin{itemize}
\item
Note that in the case $d<e<(m-1)d$ not every colouring corresponds to a tree. Therefore, consider bipartite quiver having a subquiver that is coloured as follows
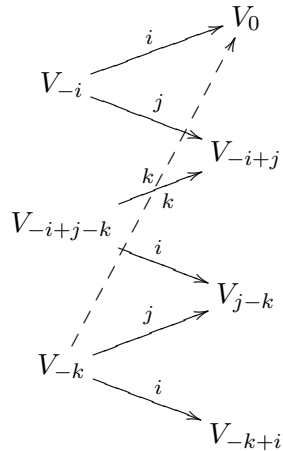
\[
\begin{xy}\xymatrix@R10pt@C30pt{
&V_0\\
V_{-i}\ar[ru]^{i}\ar[rd]^j &\\
&V_{-i+j}\\
V_{-i+j-k}\ar[ru]^k\ar[rd]^{i} &\\
&V_{j-k}\\
V_{-k}\ar@{-->}[ruuuuu]_k\ar[ru]^{j}\ar[rd]^{i} &\\
&V_{-k+i}}
\end{xy}
\]
This type of colouring induces an extra arrow so that we obviously do not get a tree. But such cycles initially break down in a second localization step, which gives us a tree module in the cases of colourings inducing such cycles. But obviously there exist colourings for which we immediately get a tree. They correspond to so called ''neighbour-avoiding-walks'' in the regular $m$-tree. In the case $m=3$ the hexagonal lattice results and $\sqrt{2}^d$ gives a lower bound. However, in general this question is unsolved.
\item In particular, this implies that there exists more than one isomorphism class of tree modules for imaginary roots. For the real roots $(0,1)$, $(1,m)$, $(m,m^2-1),\ldots$ we only get tree modules induced by the reflection of the simple representation $E_2$ considered as representations of $T(m)$.
\end{itemize}
\begin{bei}
\end{bei}
Consider the $3$-Kronecker quiver with dimension vector $(1,3)$. Then we have the simple quiver
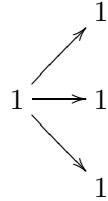
\[
\begin{xy}
\xymatrix@R20pt@C20pt{
&1\\
1\ar[ru]\ar[r]\ar[rd]&1\\&1}
\end{xy}
\]
Applying the reflection functor we get
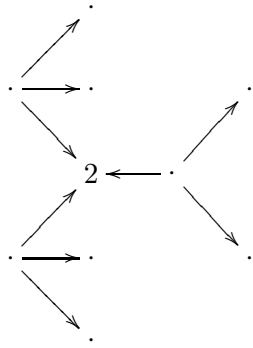
\[
\begin{xy}
\xymatrix@R20pt@C20pt{
&\cdot&&\\\cdot\ar[ru]\ar[r]\ar[rd]&\cdot&&\cdot\\&2&\cdot\ar[l]\ar[ru]\ar[rd]&\\\cdot\ar[ru]\ar[r]\ar[rd]&\cdot&&\cdot\\&\cdot&&}
\end{xy}
\]
with coefficient quiver
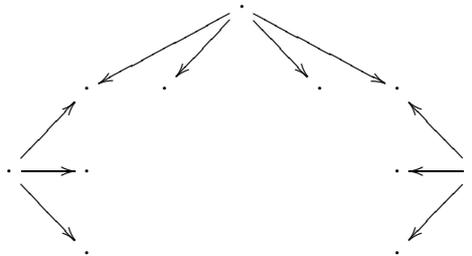
\[
\begin{xy}
\xymatrix@R20pt@C20pt{
&&&\cdot\ar[lld]\ar[ld]\ar[rd]\ar[rrd]&&&\\&\cdot&\cdot&&\cdot&\cdot&\\\cdot\ar[ru]\ar[r]\ar[rd]&\cdot&&&&\cdot&\cdot\ar[lu]\ar[l]\ar[ld]\\&\cdot&&&&\cdot&}
\end{xy}
\]
If we again apply the reflection functor, we obtain the dimension vector $(8,21)$ and we can construct the coefficient quiver as described above by considering the following quiver:
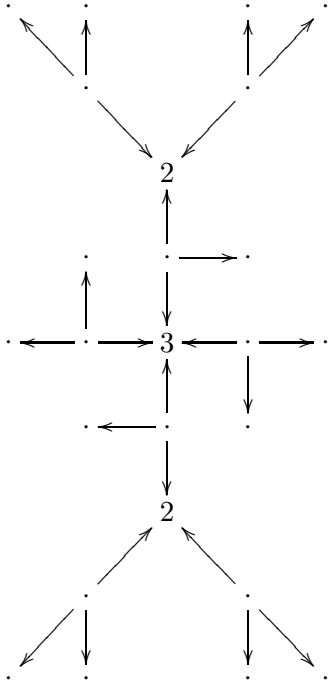
\[
\begin{xy}
\xymatrix@R20pt@C20pt{
\cdot&\cdot&&\cdot&\cdot\\
&\cdot\ar[lu]\ar[rd]\ar[u]&&\cdot\ar[u]\ar[ld]\ar[ru]\\
&&2&&\\
&\cdot&\cdot\ar[u]\ar[d]\ar[r]&\cdot\\
\cdot&\cdot\ar[u]\ar[l]\ar[r]&3&\cdot\ar[l]\ar[r]\ar[d]&\cdot\\
&\cdot&\cdot\ar[u]\ar[l]\ar[d]&\cdot\\
&&2&&\\
&\cdot\ar[ru]\ar[ld]\ar[d]&&\cdot\ar[d]\ar[rd]\ar[lu]\\
\cdot&\cdot&&\cdot&\cdot\\
}\end{xy}
\]
\begin{bei}
\end{bei}
Consider the stable simple quiver of dimension $(2,5)$ and $m=3$. Then we get
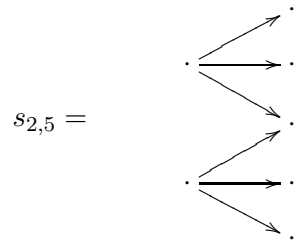
\[
\begin{xy}\xymatrix@R10pt@C30pt{
&&\cdot\\
&\cdot\ar[r]\ar[ru]\ar[rd] &\cdot\\
s_{2,5}=&&\cdot\\
&\cdot\ar[ru]\ar[rd]\ar[r]&\cdot\\
&&\cdot\\
}
\end{xy}
\]
From this we obtain the cases $(2,3)$ and $(2,4)$ by removing arrows. By applying the reflection functor, we get the dimension vector $(5,13)$:
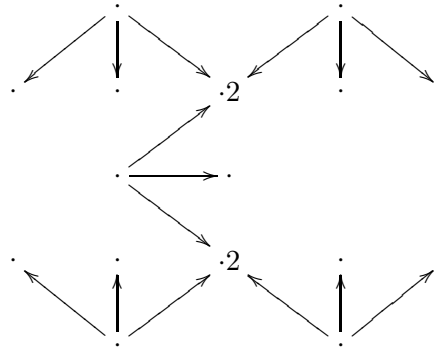
\[
\begin{xy}\xymatrix@R20pt@C30pt{
&\cdot\ar[rd]\ar[d]\ar[ld]&&\cdot\ar[rd]\ar[d]\ar[ld]&\\
\cdot&\cdot&\cdot 2&\cdot&\cdot\\
&\cdot\ar[r]\ar[ru]\ar[rd]&\cdot&&\\
\cdot&\cdot&\cdot 2&\cdot&\cdot\\
&\cdot\ar[ru]\ar[u]\ar[lu]&&\cdot\ar[ru]\ar[u]\ar[lu]&\\
}
\end{xy}
\]
From this we obtain the coefficient quiver in the way shown above. The reflected factor representations correspond to factor representations with dimension vectors $(4,10)$ and $(3,7)$ resp., thus we have:
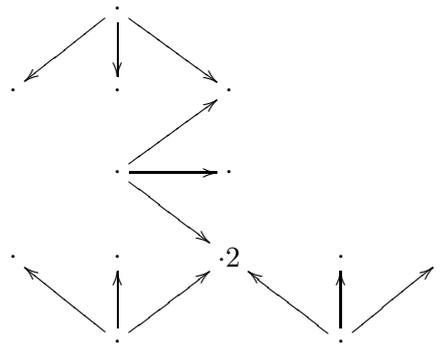
\[
\begin{xy}\xymatrix@R20pt@C30pt{
&\cdot\ar[rd]\ar[d]\ar[ld]\\
\cdot&\cdot&\cdot&&\\
&\cdot\ar[r]\ar[ru]\ar[rd]&\cdot&&\\
\cdot&\cdot&\cdot 2&\cdot&\cdot\\
&\cdot\ar[ru]\ar[u]\ar[lu]&&\cdot\ar[ru]\ar[u]\ar[lu]&\\
}
\end{xy}
\]
and
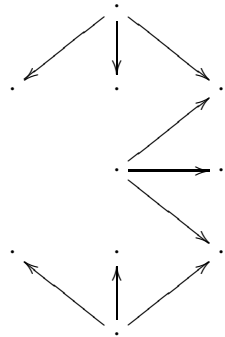
\[
\begin{xy}\xymatrix@R20pt@C30pt{
&\cdot\ar[rd]\ar[d]\ar[ld]\\
\cdot&\cdot&\cdot&&\\
&\cdot\ar[r]\ar[ru]\ar[rd]&\cdot&&\\
\cdot&\cdot&\cdot&&\\
&\cdot\ar[ru]\ar[u]\ar[lu]&&\\
}
\end{xy}
\]


\newpage
\begin{thebibliography}{Literaturverzeichnis}
\bibitem[1]{bgp}
Bernstein, Gelfand, Ponomarev: Coxeter functors and Gabriel's theorem. Russian Math. Surveys {\bf 28}, 17-32 (1973).
\bibitem[2]{fahr}
Fahr, Philipp: Infinite Gabriel-Roiter measures for the 3-Kronecker quiver. Doktorarbeit, Universität Bielefeld 2007.
\bibitem[3]{kac}
Kac, V.G.: Infinite root systems, representations of graphs and invariant theory. Inventiones mathematicae {\bf 56}, 57-92 (1980).
\bibitem[4]{kin}
King, A.: Moduli of representations of finite-dimensional algebras. Quart. J. Math. Oxford Ser. {\bf 45}, 515-530 (1994).
\bibitem[5]{rei}
Reineke, Markus: Localization in quiver moduli. Prepint 2005. To appear in: Journal für die reine und angewandte Mathematik. math.AG/0509361
\bibitem[6]{rei3}
Reineke, M.: The use of geometric and quantum group techniques for
wild quivers. {\it Representations of finite dimensional algebras
and related topics in Lie theory and geometry}, Fields
Institut Community {\bf 40}, American Math. Society, Providence, RI,
365-390 (2004).
\bibitem[7]{rei4}
Reineke, M.: Moduli of Representations of Quivers. Preprint 2008. To appear in: Proceedings of ICRA XII. arXiv:0802.2147.
\bibitem[8]{rin1}
Ringel, Claus Michael: Exceptional modules are tree modules. Linear alegbra and its Applications {\bf 275-276}, 471-493 (1998).
\bibitem[9]{rin2}
Ringel, Claus Michael: Combinatorial Representation Theory. History and future. Representations of Algebras {\bf 1}. Proceedings of the Conference on ICRA IX, Bejing 2000, Bejing Normal University Press, 122-144 (2002).
\bibitem[10]{wei}
Weist, Th.: Lokalisierung in Modulräumen. PhD thesis in preparation.
\end{thebibliography}
\end{document}